# SECOND ORDER, LINEAR AND UNCONDITIONALLY ENERGY STABLE SCHEMES FOR A HYDRODYNAMIC MODEL OF SMECTIC-A LIQUID CRYSTALS

RUI CHEN † AND XIAOFENG YANG ‡⋆ AND HUI ZHANG §

Abstract. In this paper, we consider the numerical approximations for a hydrodynamical model of smectic-A liquid crystals. The model, derived from the variational approach of the modified Oseen-Frank energy, is a highly nonlinear system that couples the incompressible Navier-Stokes equations and a constitutive equation for the layer variable. We develop two linear, second-order time-marching schemes based on the "Invariant Energy Quadratization" method for nonlinear terms in the constitutive equation, the projection method for the Navier-Stokes equations, and some subtle implicit-explicit treatments for the convective and stress terms. Moreover, we prove the well-posedness of the linear system and their unconditionally energy stabilities rigorously. Various numerical experiments are presented to demonstrate the stability and the accuracy of the numerical schemes in simulating the dynamics under shear flow and the magnetic field.

## 1. Introduction

Liquid crystals (LCs) are one important intermediate phase which exhibits features from both the solid and the fluid state, e.g., it flows like liquids, while at the same time, displays an ordering property like solid. Thus it is often viewed as the fourth state of the matter besides the gas, liquid and solid. There are two main different phases in thermotropic liquid crystals: nematic and smectic. In nematic phases, the rod-like molecules self-align to have a long-range directional order with their long axes roughly parallel. While maintaining long-range directional order, the molecules are free to flow and their center of mass positions are randomly distributed as in a liquid, see [6,7,13,17,21–24,43,50,66,71,74–77,84]. In smectic phases, the molecules maintain the general orientational order of nematics, but also tend to align themselves in layers or planes. Hence, molecules in this state show a degree of translational order that are not present in the nematic phase. Motion is restricted to within these planes, and separate planes are observed to flow past each other, see [5,7,10,11,14,24–26,30,39,43,48]. Note there are many different smectic phases, all characterized by different types and degrees of positional and orientational order. Here we consider the numerical approximations for the smectic-A phase, where the directions of molecules are perpendicular to the smectic plane, and there is no particular positional order inside the layer.

There is a large quantity of studies on modeling and simulation to investigate the flows of liquid crystal systems. One of the most well-known continuum theories is the Ericksen-Leslie theory [16, 17, 38], where, to describe this anisotropic structure, a dimensionless unit vector, called the *director*, is introduced to represent the direction of preferred orientation of molecules in the neighborhood of any point. The corresponded mathematical model can often be obtained from the energetic variational approach for the phenomenon-logical Oseen-Frank energy, leading to a well-posed nonlinear gradient flow system. The first mathematical model for the smectic-A phase is developed in the pioneering work [13] by de Gennes et. al., where the Oseen-Frank energy is modified by coupling the director field and a complex order parameter that represents the layer structure. Following the model of de Genes et. al., a number of models for smectic phase have been developed and studied during the last two decades, see [1, 7, 10, 14, 24, 30, 43, 52]. In this paper, we consider the numerical approximations for solving a particular hydrodynamics coupled smectic-A model developed by W. E. in [14] since it appears to the minimal model of unknowns, where the director field is assumed to be strictly equal to the gradient of the layer and thus the total free energy is reduced to a simplified version with one order parameter. In addition, rather than imposing non-convex constraint directly on the gradient of the layer variable, we use a commonly used technique in liquid crystal theory to modify the free energy by adding a penalization potential of a Ginzburg-Landau type. Such a term can efficiently relax the unit norm constraint numerically, while, in the meantime, it also introduces a stiffness issue into the system [2, 11, 30, 73, 84], for which certain numerical methods like fully implicit or explicit type methods (cf. [20, 55]), are numerically unstable.

† School of Science, Beijing University of Posts and Telecommunications, Beijing 100876, P.R. China, Email: ruichenbnu@gmail.com.

‡⋆Corresponding author, Department of Mathematics, University of South Carolina, Columbia, SC, 29208, USA; Beijing Institute for Scientific and Engineering Computing, Beijing University of Technology, Beijing 100124, China. Email: xfyang@math.sc.edu.

§School of Mathematical Sciences, Beijing Normal University, Laboratory of Mathematics and Complex Systems, Ministry of Education, Beijing 100875, P.R. China, Email: hzhang@bnu.edu.cn.

From the numerical point of view, for a stiff PDE system, we expect to establish schemes that can verify the so called "energy stable" property at the discrete level irrespectively of the coarseness of the discretization, namely, the energy stability does not impose any limitations on the time step. In what follows, those algorithms will be called *unconditionally energy stable.* Schemes with this property is specially preferred since it is not only critical for the numerical scheme to capture the correct long time dynamics with large time steps, but also provides sufficient flexibility for dealing with the stiffness issue. However, it is remarkable that, unlike the enormous algorithm developments on the nematic models (cf. [12, 37, 42, 45–47, 57, 85, 87–90]), very few attempts of developing energy stable schemes had been made for smectic models in any form. We notice, for solving the particular smectic-A model of W. E. [14], Guillén et. al. developed a linear, second order scheme in [30], where the nonlinear term induced by the penalization Ginzburg-Landau potential is approximated by a Hermite quadrature formula. This scheme can be regarded as one of the limited efforts in the algorithm designs, however, it is not unconditionally energy stable, i.e., there exists a time step constraint that is dependent on the penalization parameter, so it is not efficient in practice.

Therefore, in this paper, the main purpose is to develop some more efficient and effective numerical schemes for solving the hydrodynamics coupled smectic-A model in [14]. We expect that our developed schemes can own the following three desired properties, i.e., (i) accurate (second order in time); (ii) stable (the unconditional energy dissipation law holds); and (iii) easy to implement and efficient (only need to solve some fully linear equations at each time step). To achieve such a goal, instead of using traditional discretization approaches like simple implicit [20], stabilized explicit [8, 45, 54, 55, 57–59, 68, 72, 83, 85], convex splitting [67, 88, 89], or other various tricky Taylor expansions [30, 63] to discretize the nonlinear potentials, we adopt the so-called "Invariant Energy Quadratization" (IEQ) method, which is a novel approach and had been successfully applied for various *gradient flow models* in the authors' recent work (cf. [9, 69, 70, 78–82, 86, 87]). The essential idea of the IEQ method is to transform the free energy into a quadratic form (since the nonlinear potential is usually bounded from below) of a set of new variables via a change of variables. The new, equivalent system still retains the equivalent energy dissipation law in terms of the new variables. Through such a reformulation, all nonlinear terms can be treated semi-explicitly that leads to a well-posed linear system at each time step.

When this IEQ method is applied to the flow coupled model such as the smectic-A model considered in this paper, there are still new challenges due to the nonlinear couplings between the multiple variables, namely, the velocity, the director field as well as the layer, where, in particular, the equation for the velocity is not gradient flow model. To this end, we use the projection method to solve Navier-Stokes equations, and a subtle implicit-explicit treatment to treat the convective and stress terms. Finally we obtain two efficient schemes that are accurate (second order in time), easy-to-implement (linear), and unconditionally energy stable (with a discrete energy dissipation law). Moreover, we rigorously prove that the well-posedness and unconditionally energy stabilities hold for the two proposed schemes, and demonstrate the stability and the accuracy of the proposed schemes through a number of classical benchmark simulation, in particular, the layer motions under shear flow and magnetic force. To the best of the authors' knowledge, the proposed schemes here are the first second order accurate schemes for the flow coupled smectic-A model with unconditional energy stabilities.

The rest of the paper is organized as follows. In Section 2, we present the whole system and show the energy law in the continuous level. In Section 3, we develop the numerical schemes and prove their well-posedness and unconditional stabilities. In Section 4, we present various 2D numerical experiments to demonstrate the stability and the accuracy of the developed numerical schemes in simulating the dynamics under shear flow and the magnetic field. Finally, some concluding remarks are presented in Section 5.

## 2. Model

We now give a brief introduction for the hydrodynamical smectic-A phase model in [14, 30]. Let $\Omega \subset \mathbb{R}^d$ with $d = 2, 3$ be the bounded domain occupied by the LCs with boundary $\partial\Omega$. The standard Oseen-Frank distortional energy for the bulk free energy takes the following form:

$$E(\boldsymbol{d}) = \int_\Omega \Big(\frac{K_1}{2}(\nabla\cdot\boldsymbol{d})^2 + \frac{K_2}{2}(\boldsymbol{d}\cdot(\nabla\times\boldsymbol{d}))^2 + \frac{K_3}{2}|\boldsymbol{d}\times(\nabla\times\boldsymbol{d})|^2\Big)d\boldsymbol{x}, \tag{2.1}$$

where the unit vector $\boldsymbol{d}$ represents the average orientation of liquid crystal molecules and $K_1, K_2, K_3$ are elastic constants for the three canonical distortional modes: splay, twist and bending, respectively. For simplicity, we suppress the anisotropic distortional elastic modes by assuming $K_1 = K_2 = K_3 = K$. Then, the Oseen-Frank energy density reduces to the Dirichlet functional

$$E(\boldsymbol{d}) = K\int_\Omega \frac{1}{2}|\nabla\boldsymbol{d}|^2 d\boldsymbol{x}. \tag{2.2}$$

For uniaxial smectic LCs, the molecules are aligned in layers with the normal vector $\boldsymbol{n}$. More specific, for smectic-A phase, $\boldsymbol{d}$ is strictly perpendicular to the layers thus $\boldsymbol{d} = \boldsymbol{n}$. Due to the incompressibility of the layers, we have $\nabla\cdot\boldsymbol{n} = 0$, then there exists a layer function $\phi(\boldsymbol{x}, t)$ such that $\nabla\phi = \boldsymbol{n}$. In turn, the Dirichlet functional is reduced to

$$E(\phi) = K\int_\Omega \frac{1}{2}(\Delta\phi)^2 d\boldsymbol{x}, \text{ with } |\nabla\phi| = 1. \tag{2.3}$$

The norm 1 constraint applied to $|\nabla\phi|$ can bring up some additional numerical challenges in algorithm designs. A common technique to overcome it is to introduce a penalty term of the Ginzburg-Landau type $F(\nabla\phi) = \frac{1}{4\epsilon^2}(|\nabla\phi|^2 - 1)^2$ with $\epsilon \ll 1$ to regularize the distortional energy in the cores of topological defects [30, 41, 45, 87, 89, 90], where $\epsilon$ is a penalization parameter that is proportional to the size of the defect core (or zone). This regularization allows the free energy to be finite at the defect core, extending the classical Ericksen-Leslie model to handle liquid crystal flows where defects are created and annihilated in time and space. Then, the regularized elastic bulk energy density is given by

$$E(\phi) = K \int_\Omega \Big(\frac{1}{2}|\Delta\phi|^2 + \frac{(|\nabla\phi|^2-1)^2}{4\epsilon^2}\Big) d\boldsymbol{x}. \tag{2.4}$$

Assuming $\boldsymbol{u}$ is the fluid velocity field and applying the generalized Fick's law that the mass flux is proportional to the gradient of the chemical potential [3,4,44], we have the following hydrodynamical model for the smectic-A phase LCs system [14, 30]:

$$\phi_t + \nabla\cdot(\boldsymbol{u}\phi) = -Mw, \tag{2.5}$$

$$w = \frac{\delta E}{\delta \phi} = K(\Delta^2\phi - \frac{1}{\epsilon^2}\nabla\cdot(|\nabla\phi|^2-1)\nabla\phi)), \tag{2.6}$$

$$\boldsymbol{u}_t + (\boldsymbol{u}\cdot\nabla)\boldsymbol{u} - \nabla\cdot\sigma(\boldsymbol{u},\phi) + \nabla p + \phi\nabla w = 0, \tag{2.7}$$

$$\nabla\cdot\boldsymbol{u} = 0, \tag{2.8}$$

where $p$ is the pressure, $M$ is the elastic relaxation time, $\sigma$ is the dissipative (symmetric) stress tensor given in [14] that reads as,

$$\begin{aligned}\sigma(\boldsymbol{u},\phi) =& \mu_1\big(\nabla\phi^T D(\boldsymbol{u})\nabla\phi\big)\nabla\phi\otimes\nabla\phi + \mu_4 D(\boldsymbol{u})\\ &+ \mu_5\big(D(\boldsymbol{u})\nabla\phi\otimes\nabla\phi + \nabla\phi\otimes D(\boldsymbol{u})\nabla\phi\big),\end{aligned} \tag{2.9}$$

where $\mu_1, \mu_4, \mu_5$ are nonnegative parameters, and $D(\boldsymbol{u}) = \frac{1}{2}(\nabla\boldsymbol{u} + \nabla\boldsymbol{u}^T)$ is a deformation tensor. We set the non-slip boundary condition for $\boldsymbol{u}$ and the following boundary conditions for $\phi$ to remove all boundary integral terms,

$$\boldsymbol{u}|_{\partial\Omega} = 0, \quad \partial_{\boldsymbol{m}}(\Delta\phi)|_{\partial\Omega} = 0, \quad \partial_{\boldsymbol{m}}\phi|_{\partial\Omega} = 0, \tag{2.10}$$

where $\boldsymbol{m}$ is the outward normal on the boundary. It is easy to see that the equation (2.5) is mass-conserved for the layer function $\phi$, i.e., $\frac{d}{dt}\int_\Omega \phi dx = 0$.

We can easily derive the PDE energy dissipation law for the above model. Here and after, for any function $f, g \in L^2(\Omega)$, we use $(f, g) = \int_\Omega f(\boldsymbol{x})g(\boldsymbol{x})d\boldsymbol{x}$ to denote the $L^2$ inner product between functions $f(\boldsymbol{x})$ and $g(\boldsymbol{x})$, and $\|f\|^2 = (f, f)$.

By taking the $L^2$ inner product of (2.5) with $w$, of (2.6) with $\phi_t$, of (2.7) with $\boldsymbol{u}$, and summing up the obtained equalities, we can obtain

$$\begin{aligned}&\frac{d}{dt}\int_\Omega \Big(\frac{1}{2}|\boldsymbol{u}|^2 + K\big(\frac{1}{2}|\Delta\phi|^2 + \frac{(|\nabla\phi|^2-1)^2}{4\epsilon^2}\big)\Big)d\boldsymbol{x}\\ &= -\int_\Omega \Big(\mu_1(\nabla\phi^T D(\boldsymbol{u})\nabla\phi)^2 + \mu_4|D(\boldsymbol{u})|^2 + 2\mu_5|D(\boldsymbol{u})\nabla\phi|^2 + M|w|^2\Big)d\boldsymbol{x} \le 0,\end{aligned} \tag{2.11}$$

Even though the above PDE energy law is straightforward, the variable $w$ involves the fourth order derivative of $\Delta^2\phi$, and it is not convenient to use them as test functions in numerical approximations. This makes it difficult to prove the energy dissipation law in the discrete level. To overcome it, we can reformulate the momentum equation (2.7) to an equivalent form which is more applicable for numerical approximations.

Define $\dot{\phi} = \phi_t + \nabla \cdot (\boldsymbol{u}\phi), \psi = -\Delta\phi$, and notice that $w = -\frac{\dot{\phi}}{M}$, then (2.5)-(2.7) can be rewritten as,

$$\frac{\dot{\phi}}{M} = K\Delta\psi + \frac{K}{\epsilon^2}\nabla \cdot (|\nabla\phi|^2 - 1)\nabla\phi)), \tag{2.12}$$

$$\psi = -\Delta\phi, \tag{2.13}$$

$$\boldsymbol{u}_t + (\boldsymbol{u} \cdot \nabla)\boldsymbol{u} - \nabla \cdot \sigma(\boldsymbol{u}, \phi) + \nabla p - \frac{1}{M}\phi\nabla\dot{\phi} = 0. \tag{2.14}$$

$$\nabla \cdot \boldsymbol{u} = 0. \tag{2.15}$$

with the boundary conditions as

$$\boldsymbol{u}|_{\partial\Omega} = 0, \quad \partial_{\boldsymbol{m}}\psi|_{\partial\Omega} = 0, \quad \partial_{\boldsymbol{m}}\phi|_{\partial\Omega} = 0, \tag{2.16}$$

This equivalent system (2.12)-(2.15) still admits the similar energy law. We take the time derivative for (2.13) to obtain

$$\psi_t = -\Delta\phi_t, \tag{2.17}$$

Thus, by taking the $L^2$ inner product of (2.12) with $\phi_t$, of (2.17) with $K\psi$, of (2.14) with $\boldsymbol{u}$, using the incompressible condition (2.15), and summing them up, one can obtain the similar energy law as follows,

$$\begin{aligned} &\frac{d}{dt}\int_\Omega \Big(\frac{1}{2}|\boldsymbol{u}|^2 + \frac{K}{2}|\psi|^2 + \frac{K}{4\epsilon^2}(|\nabla\phi|^2 - 1)^2\Big)d\boldsymbol{x} \\ &= -\int_\Omega \Big(\mu_1(\nabla\phi^T D(\boldsymbol{u})\nabla\phi)^2 + \mu_4|D(\boldsymbol{u})|^2 + 2\mu_5|D(\boldsymbol{u})\nabla\phi|^2 + \frac{1}{M}|\dot{\phi}|^2\Big)d\boldsymbol{x} \le 0, \end{aligned} \tag{2.18}$$

Note that the above derivation is suitable in a finite dimensional approximation since the test functions $\phi_t$ and $\psi$ are both in the same subspaces as $\phi$. Hence, it allows us to design numerical schemes which satisfy the energy dissipation law in the discrete level.

**Remark 2.1.** *In [13], de Genes et. al. presented a total free energy of smectic-A phase LCs that is described by the director field $\boldsymbol{d}$ and a complex order parameter $\Psi$ that represents the average direction of molecular alignment and the layer structure, respectively. The smectic order parameter is written as $\Psi(\boldsymbol{x}) = \rho(\boldsymbol{x})e^{iq\omega(\boldsymbol{x})}$, where $\omega(\boldsymbol{x})$ is the order parameter to describe the layer structure so that $\nabla\omega$ is perpendicular to the layer, and the smectic layer density $\rho(\boldsymbol{x})$ is the mass density of the layers. Thus the total free energy proposed by de Genes et. al. reads as follows,*

$$E(\Psi, \boldsymbol{d}) = \int_\Omega \Big(C|\nabla\Psi - iq\boldsymbol{d}\Psi|^2 + K|\nabla\boldsymbol{d}|^2 + \frac{g}{2}(|\Psi|^2 - \frac{r}{g})^2\Big)d\boldsymbol{x}, \tag{2.19}$$

*where the order parameters $C, k, g, r$ are all fixed positive constants. By assuming the density $\rho(x) = r/g$ and $\phi(\boldsymbol{x}) = \frac{\omega(\boldsymbol{x})}{d}$ and rescaling other parameters, one can obtain the normalized energy as (cf. [24]),*

$$E(\phi, \boldsymbol{d}) = \int_\Omega \Big(\frac{|\nabla\phi - \boldsymbol{d}|^2}{2\eta^2} + \frac{|\nabla\boldsymbol{d}|^2}{2}\Big)d\boldsymbol{x}. \tag{2.20}$$

*where $\eta$ is a constant determined by the domain size and other parameters. Therefore, the free energy (2.3) can be viewed the approximation of the de Genes' energy when $\eta \to 0$.*

## 3. Numerical schemes

We now construct time marching schemes for solving the model system (2.12)-(2.15). Our aim is to construct schemes that are not only easy-to-implement, but also unconditionally energy stable.

Here the term easy-to-implement is referred to "linear" and "decoupled" in comparison with its counter parts: "nonlinear" and "coupled". Thus we use the IEQ approach to discretize the double well potential since it is an efficient linear approach, and the projection methods for the Navier-Stokes equation [28, 29, 62] since it can decouple the calculations of the pressure from the velocity field. The key point of IEQ method is to make the nonlinear potential quadratic. More precisely, we define an auxiliary function $U$ as

$$U = |\nabla\phi|^2 - 1, \tag{3.1}$$

thus the total energy of (4.7) turns to a new form as

$$E(\phi, U) = K\int_\Omega \Big(\frac{1}{2}(\Delta\phi)^2 + \frac{1}{4\epsilon^2}U^2\Big)d\boldsymbol{x}. \tag{3.2}$$

Then we obtain an equivalent PDE system by taking the time derivative for the new variable $U$:

$$\frac{\dot{\phi}}{M} = K\Delta\psi + \frac{K}{\epsilon^2}\nabla\cdot(U\nabla\phi)\big), \tag{3.3}$$

$$\psi = -\Delta\phi, \tag{3.4}$$

$$U_t = 2\nabla\phi\cdot\nabla\phi_t, \tag{3.5}$$

$$\boldsymbol{u}_t + (\boldsymbol{u}\cdot\nabla)\boldsymbol{u} - \nabla\cdot\sigma(\boldsymbol{u},\phi) + \nabla p - \frac{1}{M}\phi\nabla\dot{\phi} = 0, \tag{3.6}$$

$$\nabla\cdot\boldsymbol{u} = 0. \tag{3.7}$$

The boundary conditions for the new system are still (2.10) since the equation (3.5) for the new variable $U$ is simply an ODE with time. The initial conditions read as

$$\boldsymbol{u}|_{(t=0)} = \boldsymbol{u}_0,\ \phi|_{(t=0)} = \phi_0,\ U|_{(t=0)} = |\nabla\phi_0|^2 - 1. \tag{3.8}$$

It is clear that the new equivalent system (3.3)-(3.7) still retains the similar energy law. By taking the $L^2$ inner product of (3.3) with $\phi_t$, taking the time derivative of (3.4) and the $L^2$ inner product with $K\psi$, of (3.5) with $\frac{K}{2\epsilon^2}U$, of (3.6) with $\boldsymbol{u}$, using the incompressible condition (3.7), and summing the obtained equalities up, one can obtain the similar energy law as follows,

$$\frac{d}{dt}E(\boldsymbol{u},\psi,U) = -\int_\Omega \Big(\mu_1(\nabla\phi^T D(\boldsymbol{u})\nabla\phi)^2 + \mu_4|D(\boldsymbol{u})|^2 + 2\mu_5|D(\boldsymbol{u})\nabla\phi|^2 + \frac{1}{M}|\dot{\phi}|^2\Big)d\boldsymbol{x} \le 0. \tag{3.9}$$

where

$$E(\boldsymbol{u},\psi,U) = \int_\Omega \Big(\frac{1}{2}|\boldsymbol{u}|^2 + \frac{K}{2}|\psi|^2 + \frac{K}{4\epsilon^2}U^2\Big)d\boldsymbol{x} \tag{3.10}$$

**Remark 3.1.** *We emphasize that the new transformed system* (3.3)-(3.7) *is exactly equivalent to the original system* (2.12)-(2.15)*, since* (2.12) *can be easily obtained by integrating* (3.5) *with respect to the time. For the time-continuous case, the potentials in the new free energy* (3.10) *are the same as the Lyapunov functional in the original free energy of* (2.11)*. We will develop unconditionally energy stable numerical schemes for time stepping of the transformed system* (3.3)-(3.7)*, and the proposed schemes should formally follow the new energy dissipation law* (3.9) *in the discrete sense, instead of the energy law for the originated system* (2.11).

3.1. **Crank-Nicolson Scheme.** Let $\delta t > 0$ denote the time step size and set $t^n = n\,\delta t$ for $0 \le n \le N$ with the ending time $T = N\,\delta t$. We first develop a second order scheme that is based on the Crank-Nicolson, that reads as follows.

**Scheme 1.** *Having computed* $\phi^n, U^n, \boldsymbol{u}^n, p^n$*, we update* $\phi^{n+1}, U^{n+1}, \boldsymbol{u}^{n+1}, p^{n+1}$ *as follows (we compute* $\phi^1, U^1, \boldsymbol{u}^1, p^1$ *by assuming* $\phi^{-1} = \phi^0, \psi^{-1} = \psi^0 = -\Delta\phi^0, U^{-1} = U^0, \boldsymbol{u}^{-1} = \boldsymbol{u}^0, p^{-1} = p^0$ *for the initial step):*

**Step 1:**

$$\frac{1}{M}\dot{\phi}^{n+1} = K\Delta\psi^{n+\frac{1}{2}} + \frac{K}{\epsilon^2}\nabla\cdot(U^{n+\frac{1}{2}}\nabla\phi^{\star,n+\frac{1}{2}}), \tag{3.11}$$

$$\psi^{n+1} = -\Delta\phi^{n+1}, \tag{3.12}$$

$$U^{n+1} - U^n = 2\nabla\phi^{\star,n+\frac{1}{2}}\cdot(\nabla\phi^{n+1} - \nabla\phi^n), \tag{3.13}$$

$$\frac{\widetilde{\boldsymbol{u}}^{n+1} - \boldsymbol{u}^n}{\delta t} + B(\boldsymbol{u}^{\star,n+\frac{1}{2}}, \widetilde{\boldsymbol{u}}^{n+\frac{1}{2}}) - \nabla\cdot\sigma(\widetilde{\boldsymbol{u}}^{n+\frac{1}{2}}, \phi^{\star,n+\frac{1}{2}}) + \nabla p^n - \frac{1}{M}\phi^{\star,n+\frac{1}{2}}\nabla\dot{\phi}^{n+1} = 0, \tag{3.14}$$

*with the boundary conditions*

$$\widetilde{\boldsymbol{u}}^{n+1}|_{\partial\Omega} = 0, \partial_{\boldsymbol{m}}\phi^{n+1} = \partial_{\boldsymbol{m}}\psi^{n+1}|_{\partial\Omega} = 0, \tag{3.15}$$

*where*

$$\begin{cases} B(\boldsymbol{u}, \boldsymbol{v}) = (\boldsymbol{u}\cdot\nabla)\boldsymbol{v} + \frac{1}{2}(\nabla\cdot\boldsymbol{u})\boldsymbol{v}, \\ \phi^{\star,n+\frac{1}{2}} = \frac{3}{2}\phi^n - \frac{1}{2}\phi^{n-1},\ \boldsymbol{u}^{\star,n+\frac{1}{2}} = \frac{3}{2}\boldsymbol{u}^n - \frac{1}{2}\boldsymbol{u}^{n-1}, \\ \psi^{n+\frac{1}{2}} = \frac{\psi^{n+1} + \psi^n}{2},\ \widetilde{\boldsymbol{u}}^{n+\frac{1}{2}} = \frac{\widetilde{\boldsymbol{u}}^{n+1} + \boldsymbol{u}^n}{2},\ U^{n+\frac{1}{2}} = \frac{U^{n+1} + U^n}{2}, \\ \dot{\phi}^{n+1} = \frac{\phi^{n+1} - \phi^n}{\delta t} + \nabla\cdot(\widetilde{\boldsymbol{u}}^{n+\frac{1}{2}}\phi^{\star,n+\frac{1}{2}}). \end{cases} \tag{3.16}$$

**Step 2:**

$$\frac{\boldsymbol{u}^{n+1} - \widetilde{\boldsymbol{u}}^{n+1}}{\delta t} + \frac{1}{2}\nabla(p^{n+1} - p^n) = 0, \tag{3.17}$$

$$\nabla\cdot\boldsymbol{u}^{n+1} = 0,\ \boldsymbol{u}^{n+1}\cdot\boldsymbol{m}|_{\partial\Omega} = 0. \tag{3.18}$$

**Remark 3.2.** *Here, for solving the Navier-Stokes equation, we use a second order pressure correction scheme [64] to decouple the computations of pressure from that of the velocity. This projection methods are analyzed in [53] where it is shown (discrete time, continuous space) that the schemes are second order accurate for velocity in* $\ell^2(0,T;L^2(\Omega))$ *but only first order accurate for pressure in* $\ell^\infty(0,T;L^2(\Omega))$*. The loss of accuracy for pressure is due to the artificial boundary condition* (3.17) *imposed on pressure [15, 28]. One can use the rotational projection type schemes to improve the order of pressure to* 3/2*. However, how to prove the energy stability for the corresponding schemes are open questions. We also remark that the Crank-Nicolson scheme with linear extrapolation is a popular time discretization for the Navier-Stokes equation. We refer to [28, 36] and references therein for analysis on this type of discretization.*

Schemes (3.11)-(3.14) is totally linear scheme since we handle the convective and stress term by compositions of implicit (Crank-Nicolson) and explicit (second order extrapolation) discretization. Apparently, the new variable $U$ brings up some extra computational cost. But actually, we do not need to calculate $U^{n+1}$ explicitly in every step. By rewriting (3.13), we obtain

$$U^{n+\frac{1}{2}} = S^n + \nabla\phi^{\star,n+\frac{1}{2}}\cdot\nabla\phi^{n+1}, \tag{3.19}$$

where $S^n = U^n - \nabla\phi^{\star,n+\frac{1}{2}}\cdot\nabla\phi^n$. Then (3.11) and (3.14) can be written as the following system with unknowns $(\phi,\boldsymbol{u})$, where $\phi^{n+1}$ and $\widetilde{\boldsymbol{u}}^{n+1}$ are its solutions,

$$\phi+\frac{\delta t}{2}\nabla\cdot(\boldsymbol{u}\phi^{\star,n+\frac{1}{2}})+\frac{KM\delta t}{2}\Delta^2\phi-\frac{KM\delta t}{\epsilon^2}\nabla\cdot(\nabla\phi^{\star,n+\frac{1}{2}}\cdot\nabla\phi)\nabla\phi^{\star,n+\frac{1}{2}})=f_1, \tag{3.20}$$

$$\begin{aligned}\frac{\delta tM}{2}\boldsymbol{u}+\frac{\delta t^2M}{4}B(\boldsymbol{u}^{\star,n+\frac{1}{2}},\boldsymbol{u})-\frac{M\delta t^2}{4}\nabla\cdot\sigma(\boldsymbol{u},\phi^{\star,n+\frac{1}{2}})\\ -\frac{\delta t}{2}\nabla\big(\phi+\frac{\delta t}{2}\nabla\cdot(\boldsymbol{u}\phi^{\star,n+\frac{1}{2}})\big)\phi^{\star,n+\frac{1}{2}}=f_2,\end{aligned} \tag{3.21}$$

where $f_1$ and $f_2$ are given from previous time steps that read as

$$\begin{cases} f_1=\phi^n-\frac{\delta t}{2}\nabla\cdot(\boldsymbol{u}^n\phi^{\star,n+\frac{1}{2}})+\frac{KM\delta t}{2}\psi^n+\frac{KM\delta t}{\epsilon^2}\nabla\cdot(S^n\nabla\phi^{\star,n+\frac{1}{2}}),\\ f_2=\frac{M\delta t}{2}\boldsymbol{u}^n-\frac{M\delta t^2}{4}B(\boldsymbol{u}^{\star,n+\frac{1}{2}},\boldsymbol{u}^n)+\frac{\delta t^2M}{4}\nabla\cdot\sigma(\boldsymbol{u}^n,\nabla\phi^{\star,n+\frac{1}{2}})-\frac{M\delta t^2}{2}\nabla p^n\\ \qquad -\frac{\delta t}{2}\phi^{\star,n+\frac{1}{2}}\nabla\big(\phi^n-\frac{\delta t}{2}\nabla\cdot(\boldsymbol{u}^n\phi^{\star,n+\frac{1}{2}})\big).\end{cases} \tag{3.22}$$

We first show the well-posedness of the above linear system (3.20)-(3.21) as follows.

**Theorem 3.1.** *The linear system* (3.20)-(3.21) *(or* (3.11)-(3.14)*) admits a unique solution in* $(\phi,\boldsymbol{u})\in(H^2,H^1)(\Omega)$.

*Proof.* By taking the $L^2$ inner product of (3.11) with 1, we obtain

$$\int_\Omega\phi^{n+1}d\boldsymbol{x}=\int_\Omega\phi^n d\boldsymbol{x}=\cdots=\int_\Omega\phi^0 d\boldsymbol{x}. \tag{3.23}$$

Let $v_\phi=\frac{1}{|\Omega|}\int_\Omega\phi^0d\boldsymbol{x}$, and we define $\widehat{\phi}=\phi-v_\phi$. Then $\int_\Omega\widehat{\phi}d\boldsymbol{x}=0$ and $(\widehat{\phi},\boldsymbol{u})$ is the solution of the following linear system with unknowns denoted by $(\phi,\boldsymbol{u})$,

$$\phi+\frac{\delta t}{2}\nabla\cdot(\boldsymbol{u}\phi^{\star,n+\frac{1}{2}})+\frac{KM\delta t}{2}\Delta^2\phi-\frac{KM\delta t}{\epsilon^2}\nabla\cdot(\nabla\phi^{\star,n+\frac{1}{2}}\cdot\nabla\phi)\nabla\phi^{\star,n+\frac{1}{2}})=f_1-v_\phi, \tag{3.24}$$

$$\begin{aligned}\frac{\delta tM}{2}\boldsymbol{u}+\frac{\delta t^2M}{4}B(\boldsymbol{u}^{\star,n+\frac{1}{2}},\boldsymbol{u})-\frac{M\delta t^2}{4}\nabla\cdot\sigma(\boldsymbol{u},\phi^{\star,n+\frac{1}{2}})\\ -\frac{\delta t}{2}\nabla\big(\phi+\frac{\delta t}{2}\nabla\cdot(\boldsymbol{u}\phi^{\star,n+\frac{1}{2}})\big)\phi^{\star,n+\frac{1}{2}}=f_2.\end{aligned} \tag{3.25}$$

We denote the above linear system (3.24)-(3.25) as

$$\mathbb{A}\boldsymbol{X}=\mathbb{B}, \tag{3.26}$$

with $\boldsymbol{X}=(\phi,\boldsymbol{u})^T$ and $\mathbb{B}=(f_1-v_\phi,f_2)^T$.

For any $\boldsymbol{X}_1=(\phi_1,\boldsymbol{u}_1)^T$ and $\boldsymbol{X}_2=(\phi_2,\boldsymbol{u}_2)^T$ with $\int_\Omega\phi_1d\boldsymbol{x}=\int_\Omega\phi_2d\boldsymbol{x}=0$ with the boundary conditions (3.15), we have

$$\boldsymbol{X}_1^T\mathbb{A}\boldsymbol{X}_2\le C_1(\|\phi_1\|_{H^2}+\|\boldsymbol{u}_1\|_{H^1})(\|\phi_2\|_{H^2}+\|\boldsymbol{u}_2\|_{H^1}), \tag{3.27}$$

where $C_1=C(\delta t,M,\epsilon^2,K,\boldsymbol{u}^{\star,n+\frac{1}{2}},\phi^{\star,n+\frac{1}{2}},\phi^n,\mu_1,\mu_4,\mu_5)$.

For any $\boldsymbol{X} = (\phi, \boldsymbol{u})^T$ with $\int_\Omega \phi d\boldsymbol{x} = 0$, we derive

$$\begin{aligned}\boldsymbol{X}^T \mathbb{A} \boldsymbol{X} =& \|\phi + \frac{\delta t}{2} \boldsymbol{u} \nabla \phi^{\star, n+\frac{1}{2}}\|^2 + \frac{KM\delta t}{2}\|\Delta \phi\|^2 + \frac{KM\delta t}{\epsilon^2}\|\nabla \phi^{\star,n+\frac{1}{2}}\nabla\phi\|^2 + \frac{\delta t M}{2}\|\boldsymbol{u}\|^2 \\ &+ \frac{M\delta t^2}{4}\Big(\mu_1\|(\nabla\phi^{\star,n+\frac{1}{2}})^T D(\boldsymbol{u})\nabla\phi^{\star,n+\frac{1}{2}}\|^2 + \mu_4\|D(\boldsymbol{u})\|^2 + 2\mu_5\|D(\boldsymbol{u})\nabla\phi^{\star,n+\frac{1}{2}}\|^2\Big) \\ \geq& C_2(\|\phi\|_{H^2}^2 + \|\boldsymbol{u}\|_{H^1}^2),\end{aligned} \tag{3.28}$$

where $C_2 = C(\delta t, M, \epsilon^2, K, \boldsymbol{u}^{\star,n+\frac{1}{2}}, \phi^{\star,n+\frac{1}{2}}, \phi^n, \mu_4)$. Then from the Lax-Milgram theorem, we conclude the linear system (3.26) admits a unique solution $(\phi, \boldsymbol{u}) \in (H^2, H^1)(\Omega)$.

□

The energy stability of the scheme (3.11)-(3.18) is presented as follows.

**Theorem 3.2.** *The scheme* (3.11)-(3.18) *is unconditionally energy stable satisfying the following discrete energy dissipation law,*

$$\begin{aligned}E_{tot-cn2}^{n+1} = E_{tot-cn2}^{n} - \frac{\delta t}{M}\|\dot{\phi}^{n+1}\|^2 - \delta t\Big(&\mu_1\|(\nabla\phi^{\star,n+\frac{1}{2}})^T D(\widetilde{\boldsymbol{u}}^{n+\frac{1}{2}})\nabla\phi^{\star,n+\frac{1}{2}}\|^2 \\ &+ \mu_4\|D(\widetilde{\boldsymbol{u}}^{n+\frac{1}{2}})\|^2 + \mu_5\|D(\widetilde{\boldsymbol{u}}^{n+\frac{1}{2}})\nabla\phi^{\star,n+\frac{1}{2}}\|^2\Big),\end{aligned} \tag{3.29}$$

*where*

$$E_{tot-cn2}^{n} = \frac{1}{2}\|\boldsymbol{u}^n\|^2 + \frac{K}{2}\|\psi^n\|^2 + \frac{K}{4\epsilon^2}\|U^n\|^2 + \frac{\delta t^2}{8}\|\nabla p^n\|^2. \tag{3.30}$$

*Proof.* By taking the $L^2$ inner product of (3.11) with $\frac{\phi^{n+1}-\phi^n}{\delta t}$ and using integration by parts, we obtain

$$\begin{aligned}&\frac{1}{M}\|\dot{\phi}^{n+1}\|^2 + \frac{1}{M}(\nabla\dot{\phi}^{n+1}, \widetilde{\boldsymbol{u}}^{n+\frac{1}{2}}\phi^{\star,n+\frac{1}{2}}) \\ &\quad= \frac{K}{\delta t}(\Delta\psi^{n+\frac{1}{2}}, \phi^{n+1}-\phi^n) - \frac{K}{\delta t\epsilon^2}\Big(U^{n+\frac{1}{2}}\nabla\phi^{\star,n+\frac{1}{2}}, \nabla(\phi^{n+1}-\phi^n)\Big).\end{aligned} \tag{3.31}$$

We take the subtraction between $n+1$ step and $n$ step for (3.12) to obtain

$$\psi^{n+1} - \psi^n = -\Delta(\phi^{n+1} - \phi^n). \tag{3.32}$$

By taking the $L^2$ inner product of (3.32) with $\frac{K}{\delta t}\psi^{n+\frac{1}{2}}$ and using integration by parts, we obtain

$$\begin{aligned}\frac{K}{2\delta t}(\|\psi^{n+1}\|^2 - \|\psi^n\|^2) &= -\frac{K}{\delta t}(\Delta(\phi^{n+1}-\phi^n), \psi^{n+\frac{1}{2}}) \\ &= -\frac{K}{\delta t}(\phi^{n+1}-\phi^n, \Delta\psi^{n+\frac{1}{2}}).\end{aligned} \tag{3.33}$$

By taking the $L^2$ inner product of (3.13) with $\frac{K}{2\epsilon^2\delta t}U^{n+\frac{1}{2}}$, we obtain

$$\frac{K}{4\epsilon^2\delta t}(\|U^{n+1}\|^2 - \|U^n\|^2) = \frac{K}{\epsilon^2\delta t}\Big(\nabla\phi^{\star,n+\frac{1}{2}}(\nabla\phi^{n+1} - \nabla\phi^n), U^{n+\frac{1}{2}}\Big). \tag{3.34}$$

By taking the $L^2$ inner product of (3.14) with $\widetilde{\boldsymbol{u}}^{n+\frac{1}{2}}$, we obtain

$$\begin{aligned}\frac{1}{2\delta t}(\|\widetilde{\boldsymbol{u}}^{n+1}\|^2 - \|\boldsymbol{u}^n\|^2) + \Big(\sigma(\widetilde{\boldsymbol{u}}^{n+\frac{1}{2}}, \nabla\phi^{\star,n+\frac{1}{2}}), \nabla\widetilde{\boldsymbol{u}}^{n+\frac{1}{2}}\Big) + (\nabla p^n, \widetilde{\boldsymbol{u}}^{n+\frac{1}{2}}) \\ - \frac{1}{M}\Big(\phi^{\star,n+\frac{1}{2}}\nabla\dot{\phi}^{n+1}, \widetilde{\boldsymbol{u}}^{n+\frac{1}{2}}\Big) = 0.\end{aligned} \tag{3.35}$$

By taking the $L^2$ inner product of (3.17) with $\boldsymbol{u}^{n+1}$ and performing integration by parts, we have

$$\frac{1}{2\delta t}(\|\boldsymbol{u}^{n+1}\|^2 - \|\widetilde{\boldsymbol{u}}^{n+1}\|^2 + \|\boldsymbol{u}^{n+1} - \widetilde{\boldsymbol{u}}^{n+1}\|^2) = 0, \tag{3.36}$$

where we use explicitly the divergence-free condition for $\boldsymbol{u}^{n+1}$ as

$$(\nabla(p^{n+1} - p^n), \boldsymbol{u}^{n+1}) = -((p^{n+1} - p^n), \nabla \cdot \boldsymbol{u}^{n+1}) = 0. \tag{3.37}$$

We rewrite the projection step (3.17) as

$$\frac{1}{\delta t}(\boldsymbol{u}^{n+1} + \boldsymbol{u}^n - 2\widetilde{\boldsymbol{u}}^{n+\frac{1}{2}}) + \frac{1}{2}\nabla(p^{n+1} - p^n) = 0. \tag{3.38}$$

By taking the inner product of the above equation with $\frac{\delta t}{2}\nabla p^n$, one arrives at

$$\frac{\delta t}{8}\Big(\|\nabla p^{n+1}\|^2 - \|\nabla p^n\|^2 - \|\nabla(p^{n+1} - p^n)\|^2\Big) = \Big(\nabla p^n, \widetilde{\boldsymbol{u}}^{n+\frac{1}{2}}\Big). \tag{3.39}$$

On the other hand, it follows directly from (3.17) that

$$\frac{\delta t}{8}\|\nabla(p^{n+1} - p^n)\|^2 = \frac{1}{2\delta t}\|\boldsymbol{u}^{n+1} - \widetilde{\boldsymbol{u}}^{n+1}\|^2. \tag{3.40}$$

Finally, by combining (3.31), (3.33), (3.34)-(3.36), (3.39) and (3.40), we obtain

$$\begin{aligned}
&\frac{1}{M}\|\dot{\phi}^{n+1}\|^2 + \frac{K}{2\delta t}(\|\psi^{n+1}\|^2 - \|\psi^n\|^2) + \frac{K}{4\epsilon^2\delta t}(\|U^{n+1}\|^2 - \|U^n\|^2) \\
&+ \frac{\delta t}{8}(\|\nabla p^{n+1}\|^2 - \|\nabla p^n\|^2) + \frac{1}{2\delta t}(\|\boldsymbol{u}^{n+1}\|^2 - \|\boldsymbol{u}^n\|^2) \\
&+ \Big(\sigma(\widetilde{\boldsymbol{u}}^{n+\frac{1}{2}}, \nabla\phi^{\star,n+\frac{1}{2}}), \nabla\widetilde{\boldsymbol{u}}^{n+\frac{1}{2}}\Big) = 0.
\end{aligned} \tag{3.41}$$

□

**Remark 3.3.** *One can formally verify that the energy law* (3.29) *is a second order approximation of the continuous energy law* (3.9) *at time level* $t^{n+\frac{1}{2}}$*.*

**Remark 3.4.** *We notice that the idea of the IEQ approach is very simple but quite different from the traditional time marching schemes. For example, it does not require the convexity as the convex splitting approach (cf. [18]) or the boundness for the second order derivative as the linear stabilization approach (cf. [55, 56, 65]). Through a simple substitution of new variables, the complicated nonlinear potentials are transformed into quadratic forms. We summarize the great advantages of this quadratic transformations as follows: (i) this quadratization method works well for various complex nonlinear terms as long as the corresponding nonlinear potentials are bounded from below; (ii) the complicated nonlinear potential is transferred to a quadratic polynomial form which is much easier to handle; (iii) the derivative of the quadratic polynomial is linear, which provides the fundamental support for linearization method; (iv) the quadratic formulation in terms of new variables can automatically maintain this property of positivity (or bounded from below) of the nonlinear potentials.*

**Remark 3.5.** *We remark that when the nonlinear potential takes the fourth order polynomial type, e.g.* $F(\psi) = (\psi^2 - 1)^2$ *where* $\psi = \phi$ *for Cahn-Hilliard equation and* $\psi = |\nabla\phi|$ *for the smectic model in this paper or the MBE model [81], this IEQ method is exactly the same as the Lagrange multiplier method in [31, 63]. But the Lagrange multiplier method will only work the fourth order polynomial type potential since its derivative* $\psi^3$ *can be decomposed into* $\lambda(\psi)\psi$ *with* $\lambda(\psi) = |\psi|^2$ *which can be viewed as a Lagrange multiplier term. However, for other type potentials, the Lagrange multiplier method is not applicable. About the application of the IEQ approach to handle other* *type of nonlinear potentials, e.g., the logarithmic Flory-Huggins potential, or anisotropic gradient entropy, etc., we refer to the authors' other work in [69, 78, 79, 81, 86, 87].*

3.2. **Adam-bashforth Scheme.** We further develop another second order version scheme based on the backward differentiation formula with the Adam-Bashforth explicit interpolation (BDF2), that reads as follows.

**Scheme 2.** *Having computed the numerical solutions of* $(\phi, U, \boldsymbol{u}, p)$ *at* $t^n$ *and* $t^{n-1}$*, we update* $\phi^{n+1}, U^{n+1}, \boldsymbol{u}^{n+1}, p^{n+1}$ *as follows:*

**Step 1:**

$$\frac{1}{M}\dot{\phi}^{n+1} = K\Delta\psi^{n+1} + \frac{K}{\epsilon^2}\nabla\cdot(U^{n+1}\nabla\phi^{*,n+1}) \tag{3.42}$$

$$\psi^{n+1} = -\Delta\phi^{n+1}, \tag{3.43}$$

$$3U^{n+1} - 4U^n + U^{n-1} = 2\nabla\phi^{*,n+1}\cdot(3\nabla\phi^{n+1} - 4\nabla\phi^n + \nabla\phi^{n-1}) \tag{3.44}$$

$$\begin{aligned}\frac{3\widetilde{\boldsymbol{u}}^{n+1} - 4\boldsymbol{u}^n + \boldsymbol{u}^{n-1}}{2\delta t} &+ B(\boldsymbol{u}^{*,n+1}, \widetilde{\boldsymbol{u}}^{n+1}) - \nabla\cdot\sigma(\widetilde{\boldsymbol{u}}^{n+1}, \phi^{*,n+1}) + \nabla p^n \\ &- \frac{1}{M}\phi^{*,n+1}\nabla\dot{\phi}^{n+1} = 0,\end{aligned} \tag{3.45}$$

*with the boundary conditions*

$$\widetilde{\boldsymbol{u}}^{n+1}|_{\partial\Omega} = 0,\ \partial_{\boldsymbol{m}}\phi^{n+1}|_{\partial\Omega} = \partial_{\boldsymbol{m}}\psi^{n+1}|_{\partial\Omega} = 0, \tag{3.46}$$

*where*

$$\begin{cases}\boldsymbol{u}^{*,n+1} = 2\boldsymbol{u}^n - \boldsymbol{u}^{n-1}, \phi^{*,n+1} = 2\phi^n - \phi^{n-1}, \\ \dot{\phi}^{n+1} = \dfrac{3\phi^{n+1} - 4\phi^n + \phi^{n-1}}{2\delta t} + \nabla\cdot(\widetilde{\boldsymbol{u}}^{n+1}\phi^{*,n+1}).\end{cases} \tag{3.47}$$

**Step 2:**

$$3\frac{\boldsymbol{u}^{n+1} - \widetilde{\boldsymbol{u}}^{n+1}}{2\delta t} + \nabla(p^{n+1} - p^n) = 0, \tag{3.48}$$

$$\nabla\cdot\boldsymbol{u}^{n+1} = 0,\ \boldsymbol{u}^{n+1}\cdot\boldsymbol{m}|_{\partial\Omega} = 0. \tag{3.49}$$

Similar to the Crank-Nicolson scheme, one can rewrite the equations (3.5) as follows:

$$U^{n+1} = Z^n + 2\nabla\phi^{*,n+1}\cdot\nabla\phi^{n+1}, \tag{3.50}$$

where $Z^n = \frac{4U^n - U^{n-1}}{3} - 2\nabla\phi^{*,n+1}\cdot\frac{4\nabla\phi^n - \nabla\phi^{n-1}}{3}$. Then $\phi^{n+1}$ and $\widetilde{\boldsymbol{u}}^{n+1}$ are the solutions for the following system with unknowns $(\phi, \boldsymbol{u})$,

$$\phi + \frac{2\delta t}{3}\nabla\cdot(\boldsymbol{u}\phi^{*,n+1}) + \frac{2KM\delta t}{3}\Delta^2\phi - \frac{4KM\delta t}{3\epsilon^2}\nabla\cdot((\nabla\phi^{*,n+1}\cdot\nabla\phi)\nabla\phi^{*,n+1}) = g_1, \tag{3.51}$$

$$\begin{aligned}&\frac{2\delta tM}{3}\boldsymbol{u} + \frac{4\delta t^2M}{9}\boldsymbol{u}^{*,n+1}\cdot\nabla\boldsymbol{u} - \frac{4M\delta t^2}{9}\nabla\cdot\sigma(\boldsymbol{u}, \phi^{*,n+1}) \\ &\qquad - \frac{2\delta t}{3}\phi^{*,n+1}\nabla\big(\phi + \frac{2\delta t}{3}\nabla\cdot(\boldsymbol{u}\phi^{*,n+1})\big) = g_2,\end{aligned} \tag{3.52}$$

where

$$\begin{cases}g_1 = \dfrac{4\phi^n - \phi^{n-1}}{3} + \dfrac{2KM\delta t}{3\epsilon^2}\nabla\cdot(Z^n\nabla\phi^{*,n+1}), \\ g_2 = \dfrac{2M\delta t}{9}(4\boldsymbol{u}^n - \boldsymbol{u}^{n-1}) - \dfrac{4M\delta t^2}{9}\nabla p^n - \dfrac{2\delta t}{9}\phi^{*,n+1}\nabla(4\phi^n - \phi^{n-1}).\end{cases} \tag{3.53}$$

**Theorem 3.3.** *The linear system* (3.42)-(3.45) *(or* (3.51)-(3.52)*) admits a unique solution in* $(\phi, \boldsymbol{u}) \in (H^2, H^1)(\Omega)$.

*Proof.* The proof of well-posedness is similar to Theorem 3.1, thus we omit the details here. □

**Theorem 3.4.** *The scheme* (3.42)-(3.49) *is unconditionally energy stable satisfying the following discrete energy dissipation law,*

$$\begin{aligned} E^{n+1}_{tot-bdf2} \le E^n_{tot-bdf2} - \frac{\delta t}{M}\|\dot\phi^{n+1}\|^2 - \delta t\Big(\mu_1\|(\nabla\phi^{*,n+1})^T D(\widetilde{\boldsymbol{u}}^{n+1})\nabla\phi^{*,n+1}\|^2 \\ + \mu_4\|D(\widetilde{\boldsymbol{u}}^{n+1})\|^2 + \mu_5\|D(\widetilde{\boldsymbol{u}}^{n+1})\nabla\phi^{*,n+1}\|^2\Big), \end{aligned} \tag{3.54}$$

*where*

$$\begin{aligned} E^{n+1}_{tot-bdf2} =& \frac{1}{2}\Big(\frac{\|\boldsymbol{u}^{n+1}\|^2}{2} + \frac{\|2\boldsymbol{u}^{n+1} - \boldsymbol{u}^n\|^2}{2}\Big) + \frac{K}{2}\Big(\frac{\|\psi^{n+1}\|^2}{2} + \frac{\|2\psi^{n+1} - \psi^n\|^2}{2}\Big) \\ &+ \frac{K}{4\epsilon^2}\Big(\frac{\|U^{n+1}\|^2}{2} + \frac{\|2U^{n+1} - U^n\|^2}{2}\Big) + \frac{\delta t^2}{3}\|\nabla p^{n+1}\|^2. \end{aligned} \tag{3.55}$$

*Proof.* By taking the $L^2$ inner product of (3.42) with $\frac{3\phi^{n+1} - 4\phi^n + \phi^{n-1}}{2\delta t}$, we obtain

$$\begin{aligned} \frac{1}{M}\|\dot\phi^{n+1}\|^2 + \frac{1}{M}\Big(\nabla\dot\phi^{n+1}, \widetilde{\boldsymbol{u}}^{n+1}\phi^{*,n+1}\Big) = \frac{K}{2\delta t}(\Delta\psi^{n+1}, 3\phi^{n+1} - 4\phi^n + \phi^{n-1}) \\ - \frac{K}{2\delta t\epsilon^2}\Big(U^{n+1}\nabla\phi^{*,n+1}, \nabla(3\phi^{n+1} - 4\phi^n + \phi^{n-1})\Big). \end{aligned} \tag{3.56}$$

We take the subtraction of (3.43) with $n$ and $n-1$ step to obtain

$$3\psi^{n+1} - 4\psi^n + \psi^{n-1} = -\Delta(3\phi^{n+1} - 4\phi^n + \phi^{n-1}) \tag{3.57}$$

By taking the $L^2$ inner product of (3.57) with $\frac{K}{2\delta t}\psi^{n+1}$, using the integration by parts and the following identity

$$2(3a - 4b + c, a) = |a|^2 - |b|^2 + |2a - b|^2 - |2b - c|^2 + |a - 2b + c|^2, \tag{3.58}$$

we obtain

$$\begin{aligned} \frac{K}{4\delta t}(\|\psi^{n+1}\|^2 - \|\psi^n\|^2 + \|2\psi^{n+1} - \psi^n\|^2 - \|2\psi^n - \psi^{n-1}\|^2 + \|\psi^{n+1} + 2\psi^n - \psi^{n-1}\|^2) \\ = -\frac{K}{2\delta t}(3\phi^{n+1} - 4\phi^n + \phi^{n-1}, \Delta\psi^{n+1}) \end{aligned} \tag{3.59}$$

By taking the $L^2$ inner product of (3.44) with $\frac{K}{4\delta t\epsilon^2}U^{n+1}$ and applying (3.58), we obtain

$$\begin{aligned} \frac{K}{8\epsilon^2\delta t}\Big(\|U^{n+1}\|^2 - \|U^n\|^2 + \|2U^{n+1} - U^n\|^2 - \|2U^n - U^{n-1}\|^2 + \|U^{n+1} - 2U^n + U^{n-1}\|^2\Big) \\ = \frac{K}{2\delta t\epsilon^2}(\nabla\phi^{*,n+1}(3\nabla\phi^{n+1} - 4\nabla\phi^n + \nabla\phi^{n-1}), U^{n+1}). \end{aligned} \tag{3.60}$$

By taking the $L^2$ inner product of (3.14) with $\widetilde{\boldsymbol{u}}^{n+1}$, we obtain

$$\begin{aligned} (\frac{3\widetilde{\boldsymbol{u}}^{n+1} - 4\boldsymbol{u}^n + \boldsymbol{u}^{n-1}}{2\delta t}, \widetilde{\boldsymbol{u}}^{n+1}) + \Big(\sigma(\widetilde{\boldsymbol{u}}^{n+1}, \nabla\phi^{*,n+1}), \nabla\widetilde{\boldsymbol{u}}^{n+1}\Big) + (\nabla p^n, \widetilde{\boldsymbol{u}}^{n+1}) \\ - \frac{1}{M}\Big(\phi^{*,n+1}\nabla\dot\phi^{n+1}, \widetilde{\boldsymbol{u}}^{n+1}\Big) = 0. \end{aligned} \tag{3.61}$$

From (3.48), for any function $\boldsymbol{v}$ with $\nabla\cdot\boldsymbol{v} = 0$, we can derive

$$(\boldsymbol{u}^{n+1}, \boldsymbol{v}) = (\widetilde{\boldsymbol{u}}^{n+1}, \boldsymbol{v}). \tag{3.62}$$

Then for the first term in (3.61), we have

$$
\begin{aligned}
&\frac{1}{2\delta t}(3\widetilde{\boldsymbol{u}}^{n+1}-4\boldsymbol{u}^n+\boldsymbol{u}^{n-1},\widetilde{\boldsymbol{u}}^{n+1})\\
&\quad=\frac{1}{2\delta t}(3\widetilde{\boldsymbol{u}}^{n+1}-3\boldsymbol{u}^{n+1},\widetilde{\boldsymbol{u}}^{n+1})+\frac{1}{2\delta t}(3\boldsymbol{u}^{n+1}-4\boldsymbol{u}^n+\boldsymbol{u}^{n+1},\widetilde{\boldsymbol{u}}^{n+1})\\
&\quad=\frac{1}{2\delta t}(3\widetilde{\boldsymbol{u}}^{n+1}-3\boldsymbol{u}^{n+1},\widetilde{\boldsymbol{u}}^{n+1})+\frac{1}{2\delta t}(3\boldsymbol{u}^{n+1}-4\boldsymbol{u}^n+\boldsymbol{u}^{n+1},\boldsymbol{u}^{n+1})\\
&\quad=\frac{1}{2\delta t}(3\widetilde{\boldsymbol{u}}^{n+1}-3\boldsymbol{u}^{n+1},\widetilde{\boldsymbol{u}}^{n+1}+\boldsymbol{u}^{n+1})+\frac{1}{2\delta t}(3\boldsymbol{u}^{n+1}-4\boldsymbol{u}^n+\boldsymbol{u}^{n-1},\boldsymbol{u}^{n+1})\\
&\quad=\frac{3}{2\delta t}(\|\widetilde{\boldsymbol{u}}^{n+1}\|^2-\|\boldsymbol{u}^{n+1}\|^2)+\frac{1}{4\delta t}\Big(\|\boldsymbol{u}^{n+1}\|^2-\|\boldsymbol{u}^n\|^2\\
&\qquad+\|2\boldsymbol{u}^{n+1}-\boldsymbol{u}^n\|^2-\|2\boldsymbol{u}^n-\boldsymbol{u}^{n-1}\|^2+\|\boldsymbol{u}^{n+1}-2\boldsymbol{u}^n+\boldsymbol{u}^{n-1}\|^2\Big).
\end{aligned}
\tag{3.63}
$$

For the projection step, we rewrite (3.17) as

$$
\frac{3}{2\delta t}\boldsymbol{u}^{n+1}+\nabla p^{n+1}=\frac{3}{2\delta t}\widetilde{\boldsymbol{u}}^{n+1}+\nabla p^n. \tag{3.64}
$$

By squaring both sides of the above equality, we obtain

$$
\frac{9}{4\delta t^2}\|\boldsymbol{u}^{n+1}\|^2+\|\nabla p^{n+1}\|^2=\frac{9}{4\delta t^2}\|\widetilde{\boldsymbol{u}}^{n+1}\|^2+\|\nabla p^n\|^2+\frac{3}{\delta t}(\widetilde{\boldsymbol{u}}^{n+1},\nabla p^n), \tag{3.65}
$$

namely, we have

$$
\frac{3}{4\delta t}(\|\boldsymbol{u}^{n+1}\|^2-\|\widetilde{\boldsymbol{u}}^{n+1}\|^2)+\frac{\delta t}{3}(\|\nabla p^{n+1}\|^2-\|\nabla p^n\|^2)=(\widetilde{\boldsymbol{u}}^{n+1},\nabla p^n). \tag{3.66}
$$

By taking the $L^2$ inner product of (3.48) with $\boldsymbol{u}^{n+1}$, we have

$$
\frac{3}{4\delta t}\Big(\|\boldsymbol{u}^{n+1}\|^2-\|\widetilde{\boldsymbol{u}}^{n+1}\|^2+\|\boldsymbol{u}^{n+1}-\widetilde{\boldsymbol{u}}^{n+1}\|^2\Big)=0. \tag{3.67}
$$

Finally, by combining (3.56), (3.57), (3.60), (3.61), (3.63), (3.66) and (3.67), we obtain

$$
\begin{aligned}
&\frac{1}{M}\|\dot{\phi}^{n+1}\|^2+\frac{3}{4\delta t}\|\boldsymbol{u}^{n+1}-\widetilde{\boldsymbol{u}}^{n+1}\|^2+\frac{\delta t}{3}(\|\nabla p^{n+1}\|^2-\|\nabla p^n\|^2)\\
&+\frac{K}{4\delta t}\Big(\|\psi^{n+1}\|^2-\|\psi^n\|^2+\|2\psi^{n+1}-\psi^n\|^2-\|2\psi^n-\psi^{n-1}\|^2+\|\psi^{n+1}-2\psi^n+\psi^{n-1}\|^2\Big)\\
&+\frac{K}{8\epsilon^2\delta t}\Big(\|U^{n+1}\|^2-\|U^n\|^2+\|2U^{n+1}-U^n\|^2-\|2U^n-U^{n-1}\|^2+\|U^{n+1}-2U^n+U^{n-1}\|^2\Big)\\
&+\frac{1}{4\delta t}\Big(\|\boldsymbol{u}^{n+1}\|^2-\|\boldsymbol{u}^n\|^2+\|2\boldsymbol{u}^{n+1}-\boldsymbol{u}^n\|^2-\|2\boldsymbol{u}^n-\boldsymbol{u}^{n-1}\|^2+\|\boldsymbol{u}^{n+1}-2\boldsymbol{u}^n+\boldsymbol{u}^{n-1}\|^2\Big)\\
&+(\sigma(\widetilde{\boldsymbol{u}}^{n+1},\nabla\phi^{*,n+1}),\nabla\widetilde{\boldsymbol{u}}^{n+1})=0,
\end{aligned}
$$

that concludes the theorem. □

**Remark 3.6.** *Heuristically, the* $\frac{1}{\delta t}(E^{n+1}_{tot-bdf2}-E^n_{tot-bdf2})$ *is a second order approximation of* $\frac{d}{dt}E(\phi,U)$ *at* $t=t^{n+1}$. *For instance, for any smooth variable* $S$ *with time, one can write*

$$
\begin{aligned}
&\Big(\frac{\|S^{n+1}\|^2+\|2S^{n+1}-S^n\|^2}{2\delta t}\Big)-\Big(\frac{\|S^n\|^2+\|2S^n-S^{n-1}\|^2}{2\delta t}\Big)\\
&\qquad\cong\Big(\frac{\|S^{n+2}\|^2-\|S^n\|^2}{2\delta t}\Big)+O(\delta t^2)\cong\frac{d}{dt}\|S(t^{n+1})\|^2+O(\delta t^2).
\end{aligned}
$$

**Remark 3.7.** *Although we consider only time discrete schemes in this paper, the results here can be carried over to any consistent finite-dimensional Galerkin type approximations since the analyses are based on the variational formulation with all test functions in the same space as the space of the trial functions. The details for the fully discrete scheme will be left to the interested readers.*

**Remark 3.8.** *For the numerical schemes proposed in this paper, the energy stability is formally derived. The error estimates for the second order scheme for the layer variable is straightforward when the velocity field is null. This is because the $H^2$ bound exists for $\phi$ from the Poincaré inequality, and the corresponding convergence analysis can be further carried out. For the hydrodynamics coupled model, we have combine the analysis work for the projection method, see [53, 64], and follow the same lines as [19] to handle the nonlinear convective and stress terms where the basic tool is to use sobolev embeddings among various Banach spaces. We will implement the rigorous error analysis in the future work.*

## 4. Numerical simulations

We now present various numerical experiments to validate the theoretical results derived in the previous sections and demonstrate the stability and accuracy of the proposed numerical schemes. In all examples, we use the inf-sup stable Iso-$P2/P1$ element [61] for the velocity and pressure, and linear element for the phase function $\phi$ and $\psi$. As for as the stable element for the Navier-Stokes variables $(u, p)$, one can read the related literatures in [32–35, 40]. If not explicit specified, the model parameters take default values given below:

$$\epsilon = 0.05, \mu_4 = 0.02, \mu_1 = \mu_5 = 0, M = 1 \times 10^{-6}, K = 0.01. \tag{4.1}$$

### 4.1. Accuracy test.

We first perform numerical simulations to test the convergence rates of the two proposed schemes (3.11)-(3.18) (denoted by CN2), and (3.42)-(3.49) (denoted by BDF2).

| $\delta t$ | $\text{Error}_u$ | Order | $\text{Error}_v$ | Order | $\text{Error}_p$ | Order | $\text{Error}_\phi$ | Order |
|---|---|---|---|---|---|---|---|---|
| $1 \times 10^{-2}$ | $4.61 \times 10^{-4}$ | – | $4.75 \times 10^{-4}$ | – | $1.18 \times 10^{-1}$ | – | $1.81 \times 10^{-4}$ | – |
| $5 \times 10^{-3}$ | $1.15 \times 10^{-4}$ | 2.003 | $1.19 \times 10^{-4}$ | 1.997 | $5.87 \times 10^{-2}$ | 1.007 | $4.53 \times 10^{-5}$ | 1.998 |
| $2.5 \times 10^{-3}$ | $2.88 \times 10^{-5}$ | 1.997 | $2.97 \times 10^{-5}$ | 2.002 | $2.94 \times 10^{-2}$ | 0.997 | $1.13 \times 10^{-5}$ | 2.003 |
| $1.25 \times 10^{-3}$ | $7.21 \times 10^{-6}$ | 1.998 | $7.43 \times 10^{-6}$ | 1.999 | $1.47 \times 10^{-2}$ | 1.000 | $2.83 \times 10^{-6}$ | 1.997 |
| $6.25 \times 10^{-4}$ | $1.80 \times 10^{-6}$ | 2.002 | $1.86 \times 10^{-6}$ | 1.998 | $7.30 \times 10^{-3}$ | 1.009 | $7.08 \times 10^{-7}$ | 1.999 |

Table 1. The $L^2$ errors for the velocity field $\boldsymbol{u} = (u, v)$, the phase variable $\phi$ and the pressure $p$ at $t = 1$ for by the scheme CN2 using different temporal resolutions with the exact solution of (4.2).

| $\delta t$ | $\text{Error}_u$ | Order | $\text{Error}_v$ | Order | $\text{Error}_p$ | Order | $\text{Error}_\phi$ | Order |
|---|---|---|---|---|---|---|---|---|
| $1 \times 10^{-2}$ | $4.00 \times 10^{-3}$ | – | $4.12 \times 10^{-3}$ | – | $3.92 \times 10^{-1}$ | – | $1.57 \times 10^{-3}$ | – |
| $5 \times 10^{-3}$ | $9.62 \times 10^{-4}$ | 2.055 | $9.91 \times 10^{-4}$ | 2.055 | $1.96 \times 10^{-1}$ | 1.000 | $3.78 \times 10^{-4}$ | 2.054 |
| $2.5 \times 10^{-3}$ | $2.36 \times 10^{-4}$ | 2.027 | $2.43 \times 10^{-4}$ | 2.027 | $9.92 \times 10^{-2}$ | 0.982 | $9.26 \times 10^{-5}$ | 2.029 |
| $1.25 \times 10^{-3}$ | $5.83 \times 10^{-5}$ | 2.017 | $6.00 \times 10^{-5}$ | 2.017 | $4.91 \times 10^{-2}$ | 1.014 | $2.29 \times 10^{-5}$ | 2.015 |
| $6.25 \times 10^{-4}$ | $1.45 \times 10^{-5}$ | 2.007 | $1.49 \times 10^{-5}$ | 2.009 | $2.46 \times 10^{-2}$ | 0.997 | $5.69 \times 10^{-6}$ | 2.008 |

TABLE 2. The $L^2$ errors for the velocity field $\boldsymbol{u} = (u, v)$, the phase variable $\phi$ and the pressure $p$ at $t = 1$ for by the scheme BDF2 using different temporal resolutions with the exact solution of (4.2).

| $\delta t$ | $\text{Error}_u$ | Order | $\text{Error}_v$ | Order | $\text{Error}_p$ | Order | $\text{Error}_\phi$ | Order |
|---|---|---|---|---|---|---|---|---|
| $1 \times 10^{-2}$ | $4.61 \times 10^{-4}$ | – | $4.75 \times 10^{-4}$ | – | $1.18 \times 10^{-1}$ | – | $2.01 \times 10^{-4}$ | – |
| $5 \times 10^{-3}$ | $1.18 \times 10^{-4}$ | 1.966 | $1.21 \times 10^{-4}$ | 1.972 | $5.87 \times 10^{-2}$ | 1.007 | $5.12 \times 10^{-5}$ | 1.973 |
| $2.5 \times 10^{-3}$ | $3.12 \times 10^{-5}$ | 1.919 | $3.21 \times 10^{-5}$ | 1.914 | $2.94 \times 10^{-2}$ | 0.997 | $1.23 \times 10^{-5}$ | 2.057 |
| $1.25 \times 10^{-3}$ | $8.14 \times 10^{-6}$ | 1.938 | $8.21 \times 10^{-6}$ | 1.967 | $1.47 \times 10^{-2}$ | 1.000 | $2.97 \times 10^{-6}$ | 2.050 |
| $6.25 \times 10^{-4}$ | $2.09 \times 10^{-6}$ | 1.961 | $2.16 \times 10^{-6}$ | 1.926 | $7.33 \times 10^{-3}$ | 1.003 | $7.92 \times 10^{-7}$ | 1.901 |

TABLE 3. The $L^2$ numerical errors at $t = 1$ that are computed by the scheme CN2 using various temporal resolutions with the initial conditions of (4.3), for mesh refinement test in time.

| $\delta t$ | $\text{Error}_u$ | Order | $\text{Error}_v$ | Order | $\text{Error}_p$ | Order | $\text{Error}_\phi$ | Order |
|---|---|---|---|---|---|---|---|---|
| $1 \times 10^{-2}$ | $4.00 \times 10^{-3}$ | – | $4.12 \times 10^{-3}$ | – | $3.93 \times 10^{-1}$ | – | $1.56 \times 10^{-3}$ | – |
| $5 \times 10^{-3}$ | $9.64 \times 10^{-4}$ | 2.052 | $9.93 \times 10^{-4}$ | 2.052 | $1.96 \times 10^{-1}$ | 1.003 | $3.71 \times 10^{-4}$ | 2.072 |
| $2.5 \times 10^{-3}$ | $2.39 \times 10^{-4}$ | 2.012 | $2.46 \times 10^{-4}$ | 2.013 | $9.82 \times 10^{-2}$ | 0.997 | $9.95 \times 10^{-5}$ | 1.898 |
| $1.25 \times 10^{-3}$ | $5.74 \times 10^{-5}$ | 2.057 | $5.94 \times 10^{-5}$ | 2.050 | $4.91 \times 10^{-2}$ | 1.000 | $2.60 \times 10^{-5}$ | 1.936 |
| $6.25 \times 10^{-4}$ | $1.49 \times 10^{-5}$ | 1.945 | $1.57 \times 10^{-5}$ | 1.919 | $2.46 \times 10^{-2}$ | 0.997 | $6.42 \times 10^{-6}$ | 2.017 |

TABLE 4. The $L^2$ numerical errors at $t = 1$ that are computed by the scheme BDF2 using various temporal resolutions with the initial conditions of (4.3), for mesh refinement test in time.

4.1.1. *Presumed exact solution.* In the first example, we set the computed domain to be $\Omega = [0, 2]^2$ and assume the following functions

$$
\begin{cases}
u(t, x, y) = \pi \sin(2\pi y) \sin^2(\pi x) \sin t, \\
v(t, x, y) = -\pi \sin(2\pi x) \sin^2(\pi y) \sin t, \\
\phi(t, x, y) = 2 + \cos(\pi x) \cos(\pi y) \sin t, \\
p(t, x, y) = \cos(\pi x) \sin(\pi y) \sin t
\end{cases}
\tag{4.2}
$$

to be the exact solution, and impose some suitable force fields such that the given solution can satisfy the system. We use 10145 nodes and 19968 triangle elements for the discretization of the space. In Table 1 and 2, we list the $L^2$ errors of the velocity field $\boldsymbol{u} = (u, v)$, the phase variable $\phi$ and the pressure $p$ between the numerically simulated solution and the exact solution at $t = 1$ with different time step sizes, for the schemes CN2 and BDF2, respectively. We observe that the schemes CN2 and BDF2 achieve almost perfect second order accuracy for $\boldsymbol{u}$ and $\phi$, and first order accuracy for $p$ in time as expected, respectively.

4.1.2. *Mesh refinement in time.* We now perform more refinement tests for temporal convergence. We set the initial conditions as follows,

$$\phi_0 = 2,\ \boldsymbol{u}_0 = (u_0, v_0) = 0,\ p_0 = 0. \tag{4.3}$$

We perform the refinement test of the time step size. Since the exact solutions are not known, we choose the solution obtained by the scheme CN2 with the time step size $\delta t = 1 \times 10^{-6}$ as the benchmark solution for computing errors. We present the $L^2$ of the variables between the numerical solution and the exact solution at $t = 1$ with different time step sizes in Table 3 and Table 4 for the schemes CN2 and BDF2, respectively. As the previous tests, we observe that the schemes CN2 and BDF2 achieve almost perfect second order accuracy for $\boldsymbol{u}$ and $\phi$, and the first order accuracy for $p$, respectively.

4.2. **Layer motion.** In this example, we consider the layer motion using the second order scheme CN2. The following initial conditions are taken as follows,

$$\phi_0 = \sin x \cos^2 y, \boldsymbol{u}_0 = (u_0, v_0) = 0, p_0 = 0, \tag{4.4}$$

that had been studied in [30]. We set the computed domain to be $\Omega = [-1, 1]^2$ and the space is discretized by using 10145 nodes and 19968 triangle elements. The model parameters are from (4.1).

We emphasize that any time step size $\delta t$ is allowable for the computations from the stability concern since all developed schemes are unconditionally energy stable. But larger time step will definitely induce large numerical errors. Therefore, we need to discover the rough range of the allowable maximum time step size in order to obtain good accuracy and to consume as low computational cost as possible. This time step range could be estimated through the energy evolution curve plots, shown in Fig. 1, where we compare the time evolution of the free energy for five different time step sizes until $t = 200$ using the second order scheme CN2. We observe that all five energy curves show decays monotonically for all time step sizes, which numerically confirms that our algorithms are unconditionally energy stable. For smaller time steps of $\delta t = 0.0001, 0.0005, 0.001, 0.005, 0.01$, all five energy curves coincide very well, that means we can just use the maximum allowable time step $\delta t = 0.01$ without worrying the accuracy.

In Fig. 2 and Fig. 3, we show the dynamical evolution of the layer function $\phi$, and the velocity field $\boldsymbol{u}$ until the simulation reaches the steady state, respectively. The obtained results show qualitatively consistent features with the numerical examples in [30].

4.3. **Layer undulation under shear flow.** In this example, we consider the numerical simulations of the layer undulation under shear flow using the second order scheme CN2. We set the computed domain to be $\Omega = [-1, 1] \times [-0.5, 0.5]$ and the space is discretized by using 10145 nodes and 19968 triangle elements. The initial condition reads as follows:

$$\phi_0 = y, \boldsymbol{u}_0 = (0.4y, 0), p_0 = 0, \tag{4.5}$$

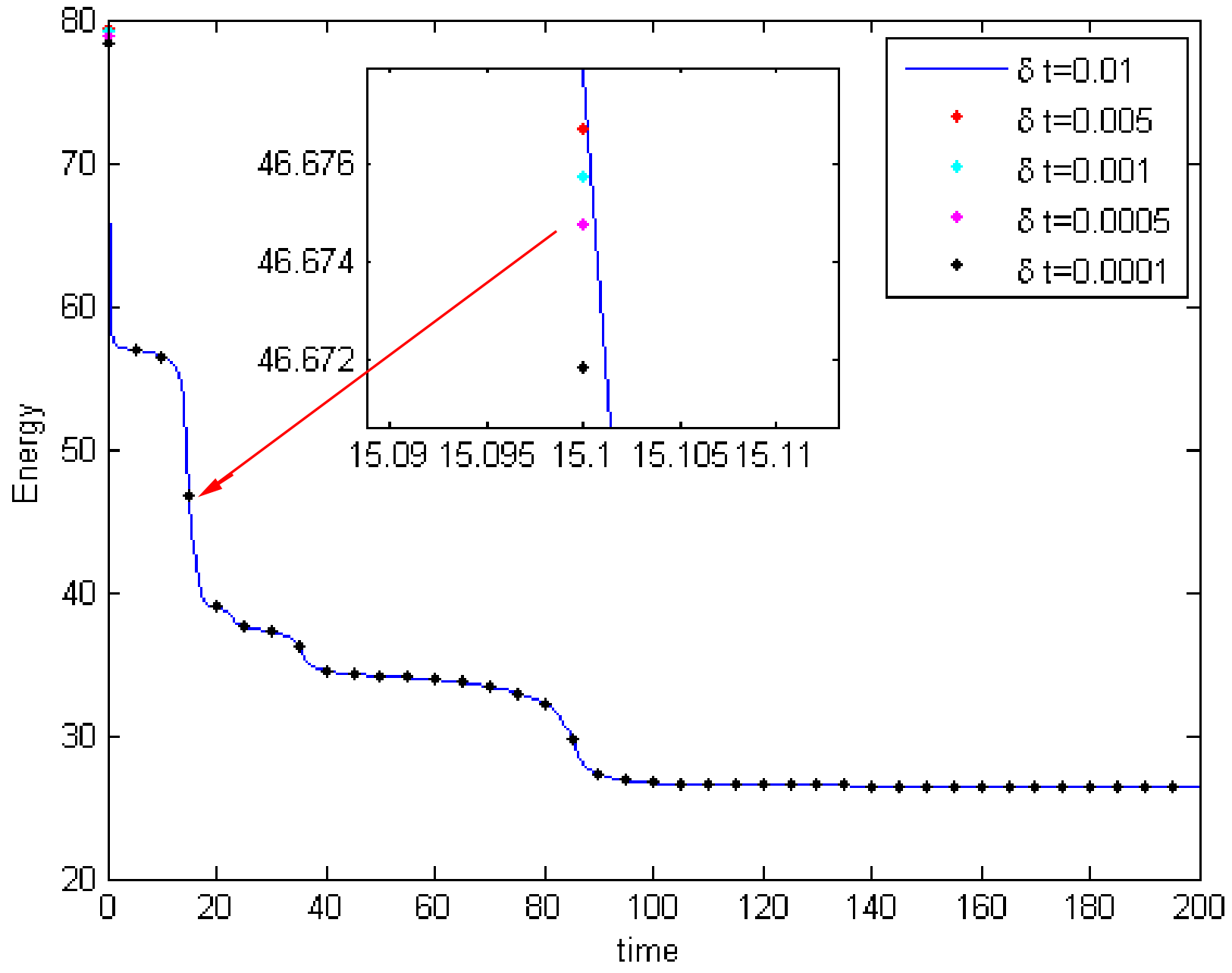

FIG. 1. Time evolution of the free energy functional till $t = 200$ for five different time steps of $\delta t = 0.01$, 0.005, 0.001, 0.0005, and 0.0001 using the scheme CN2. The energy curves show the decays for all time steps, which confirms that our algorithm is unconditionally stable. The small differences in the energy evolution for all five time steps are shown as well.

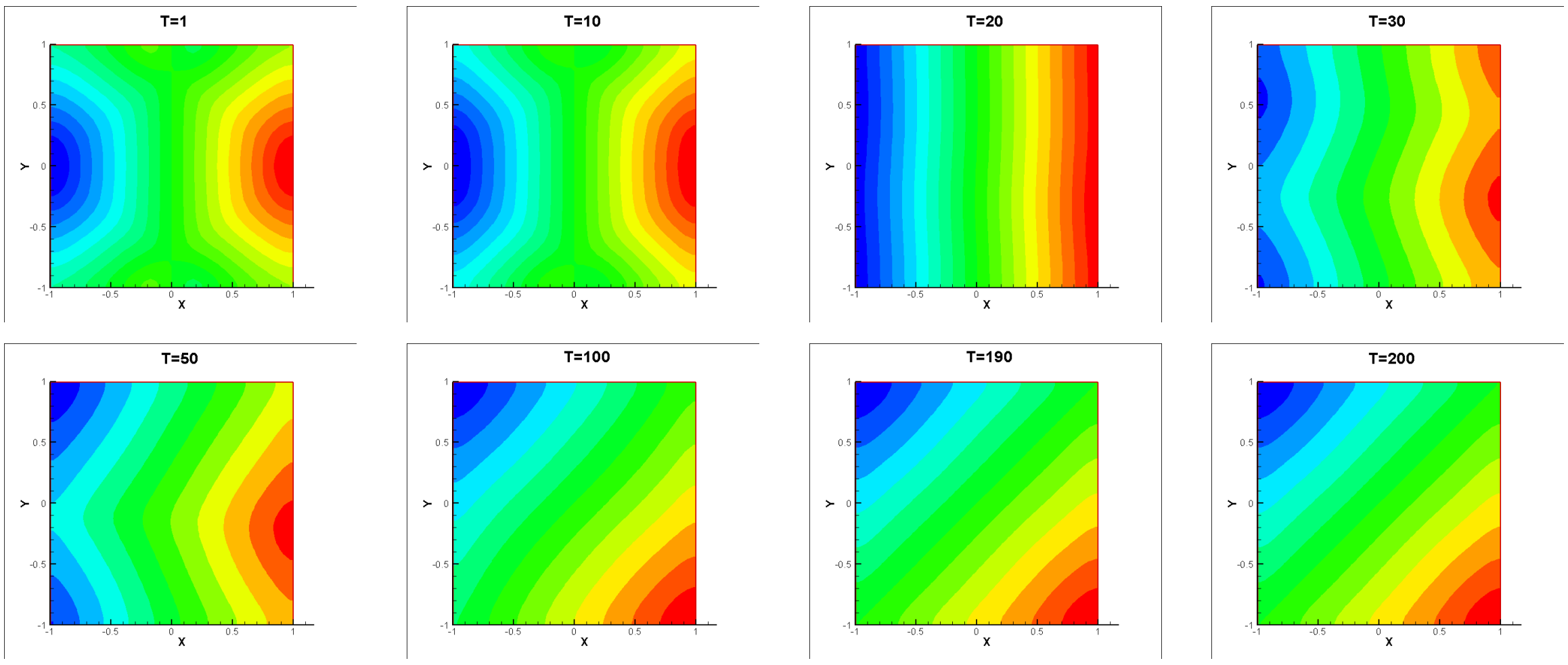

FIG. 2. The dynamical evolution of the layer function $\phi$ for the layer motion example with the time step $\delta t = 0.01$. Snapshots of the numerical approximation are taken at $t = 1$, 10, 20, 30. 50, 100, 190, and 200.

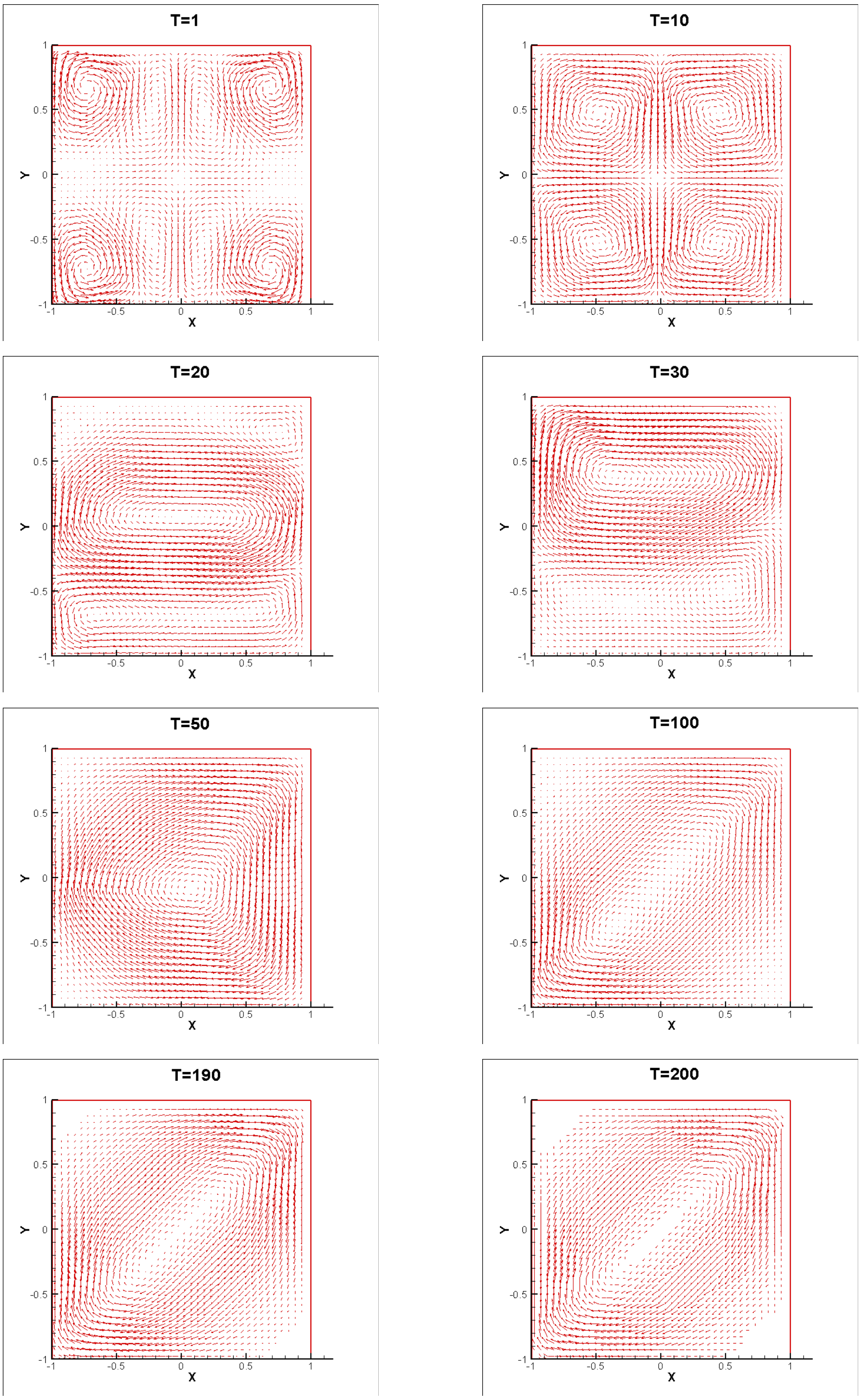

Fig. 3. The profiles of the velocity field $\boldsymbol{u}$ for the layer motion example with the time step $\delta t = 0.01$. Snapshots of the numerical approximation of $\boldsymbol{u}$ are taken at $t = 1$, 10, 20, 30. 50, 100, 190, and 200.

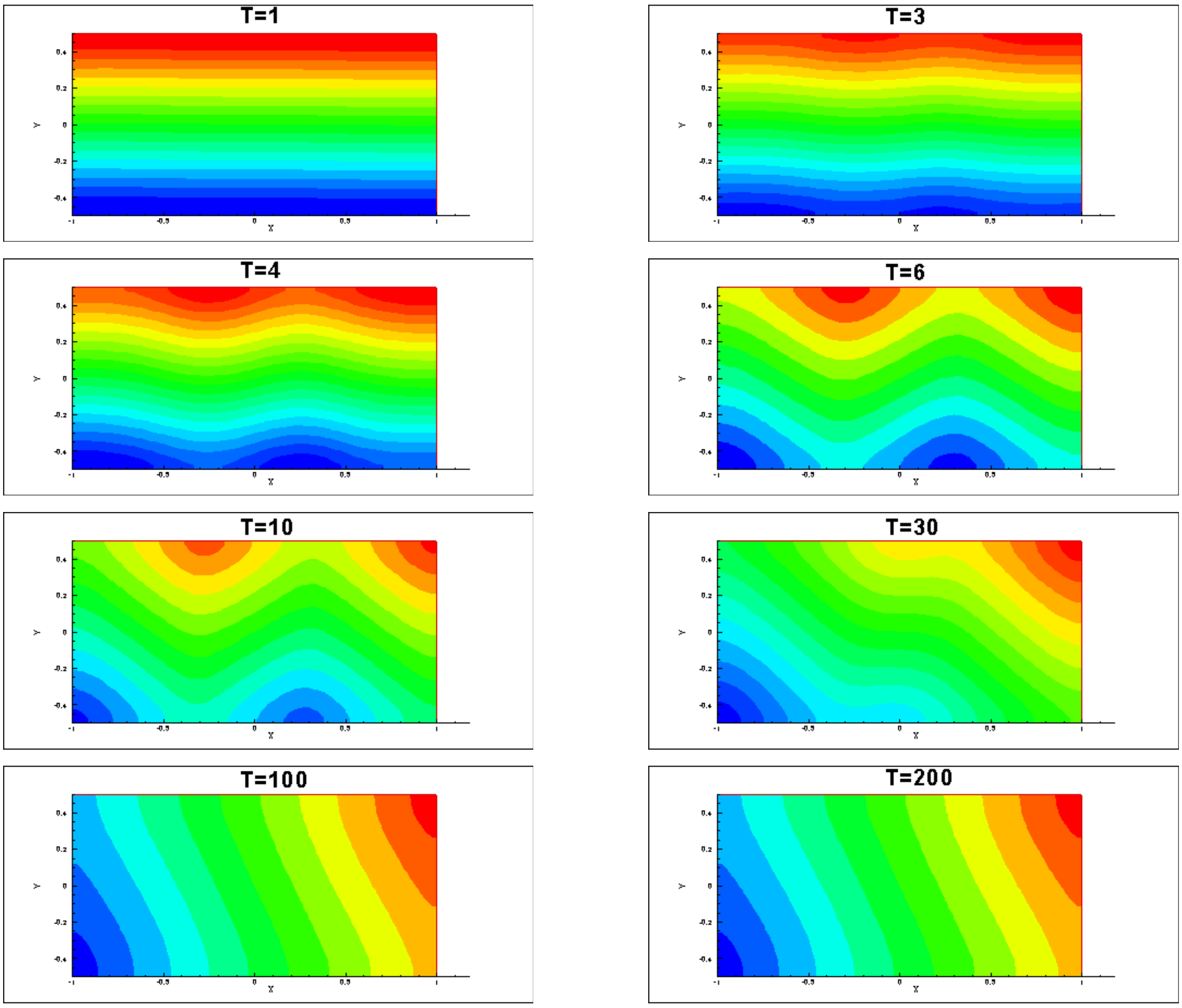

FIG. 4. The dynamical evolution of the layer function $\phi$ for the layer undulations under the shear flow with the time step $\delta t = 0.01$. Snapshots of the numerical approximation are taken at $t = 1, 3, 4, 6, 10, 30, 100, 200$.

and the boundary condition for the velocity field are set to be

$$\boldsymbol{u}|_{y=0.5} = (0.2, 0), \ \boldsymbol{u}|_{y=-0.5} = (-0.2, 0), \ \boldsymbol{u}|_{x=\pm 1} = (0, 0). \tag{4.6}$$

The model parameters are still from (4.1).

In Fig. 4 and Fig. 5, we show the dynamical evolution of the layer function $\phi$, and the velocity field $\boldsymbol{u}$ until the simulation reaches the steady state, respectively. The obtained profiles of undulational layers are consistent with the theoretical results predicted in [48] and the numerical results using the molecular dynamics approach in [60].

4.4. **The sawtooth feature under external magnetic field.** Applying an external magnetic field is one of the most efficient approach to control and produce various nano-structured materials, and had been well-studied in a number of experimental, modeling and numerical literatures, see [25,27,49,51]. In the last numerical example, we consider the dynamical behaviors of the smectic-A LCs in the presence of an applied magnetic field. When an external magnetic filed is applied, an additional term contributed by it is added to the free energy of the model system, that reads as

$$E(\phi, \boldsymbol{u}) = \int_\Omega \Big(\frac{1}{2}|\boldsymbol{u}|^2 + \frac{K}{2}|\Delta\phi|^2 + K\frac{(|\nabla\phi|^2 - 1)^2}{4\epsilon^2} - \tau(\nabla\phi \cdot \boldsymbol{h})^2\Big) d\boldsymbol{x}, \tag{4.7}$$

where $\boldsymbol{h}$ is a given unit vector representing the direction of the magnetic field, and $\tau$ is a nonnegative parameter denoting the strength of the applied magnetic field.

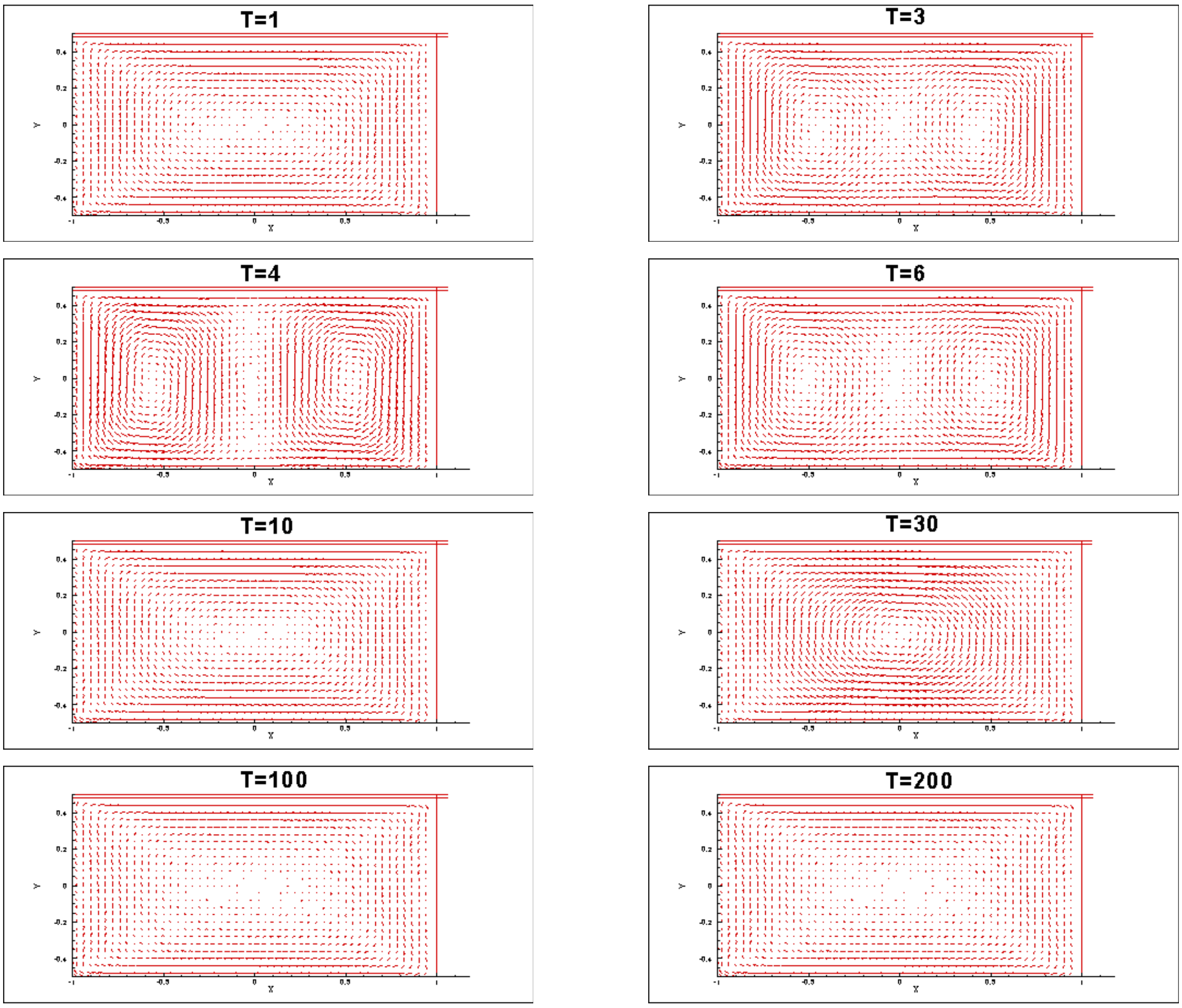

FIG. 5. The dynamical evolution of the velocity field $\boldsymbol{u}$ for the layer undulations under the shear flow with the time step $\delta t = 0.01$. Snapshots of the numerical approximation are taken at $t = 1, 3, 4, 6, 10, 30, 100, 200$.

Thus the new equation for the layer function $\phi$ reads as follows

$$\phi_t + \nabla \cdot (\boldsymbol{u}\phi) = -Mw, \tag{4.8}$$

$$w = \frac{\delta E}{\delta \phi} = K(\Delta^2\phi - \frac{1}{\epsilon^2}\nabla \cdot (|\nabla\phi|^2 - 1)\nabla\phi)) + \tau\nabla \cdot (\nabla\phi \cdot \boldsymbol{h})\boldsymbol{h}, \tag{4.9}$$

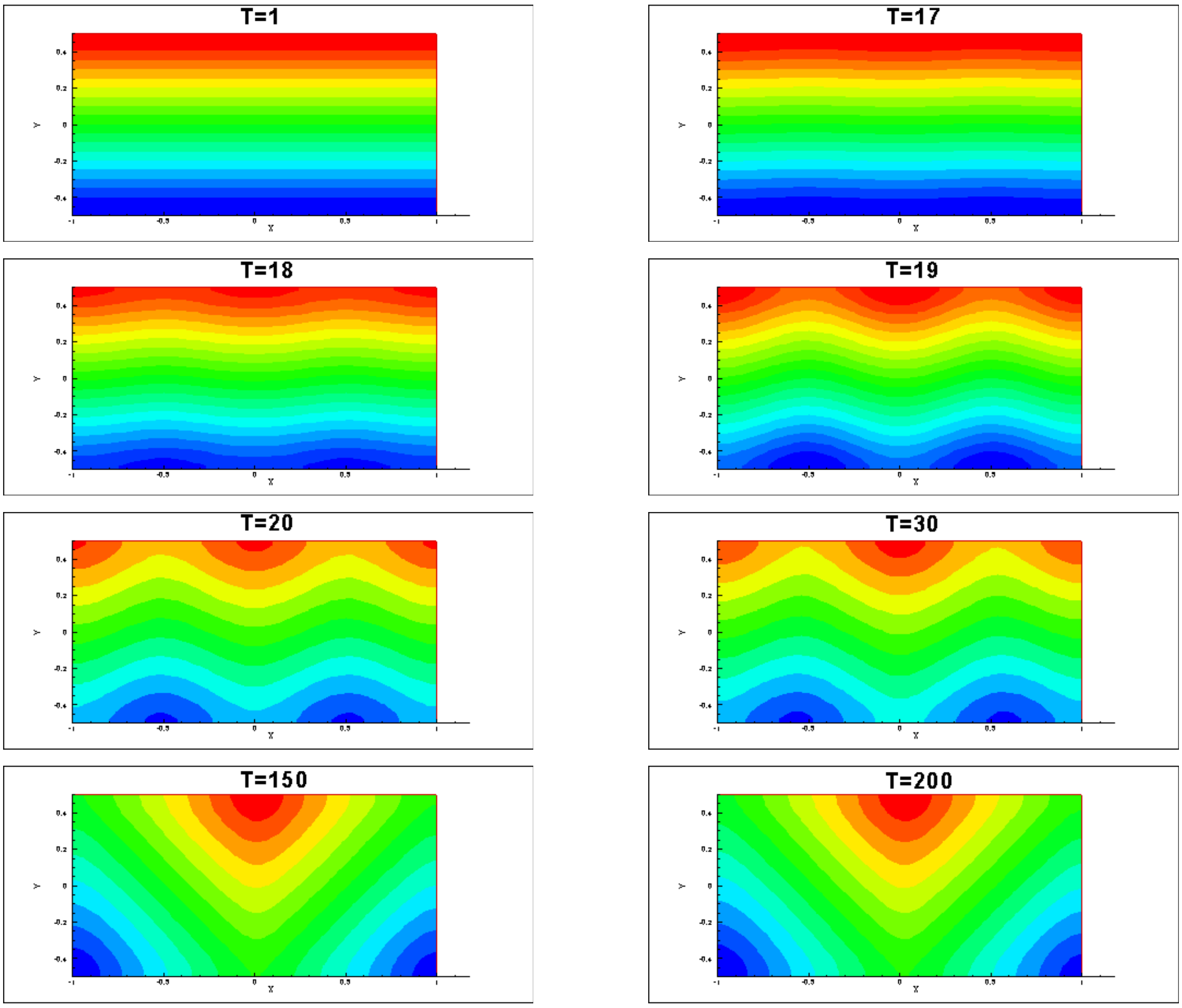

FIG. 6. The dynamical evolution of the layer function $\phi$ under external the magnetic field with the time step $\delta t = 0.01$. Snapshots of the numerical approximation are taken at $t = 1$, 17, 18, 19, 20, 30, 150, 200.

and the equations for the fluid velocity are still (2.7)-(2.8). The magnetic field term can be viewed as an imposed external force, i.e., we treat this term by the second order extrapolations.

We let $\boldsymbol{h} = (1, 0)$ and $\tau = 10$ and choose the same initial conditions, computed domain and the space discretizations as the previous shear flow example. In Fig. 6, we present that the dynamical motion of the layer variable $\phi$, the undulation profile is formed from $t = 2.1$ to $t = 3$. This sawtooth feature is qualitatively consistent with the numerical simulation in [27] using the de Gennes' smectic-A model. The final equilibrium solution is obtained after $t = 150$. We present the snapshots of the velocity field in Fig. 7 as well.

## 5. Conclusions and remarks

In this paper, we have constructed a set of efficient numerical schemes for solving the hydrodynamics coupled smectic-A LCs model. The schemes are (i) second order accurate in time; (ii) unconditional energy stable; and (iii) linear and easy to implement. Various numerical results are presented to validate the accuracy of our schemes. We have also presented a number of numerical simulations to show the morphological evolutions, in particular, the layer undulation under shear flow as well as the sawtooth profile induced by the external magnetic field.

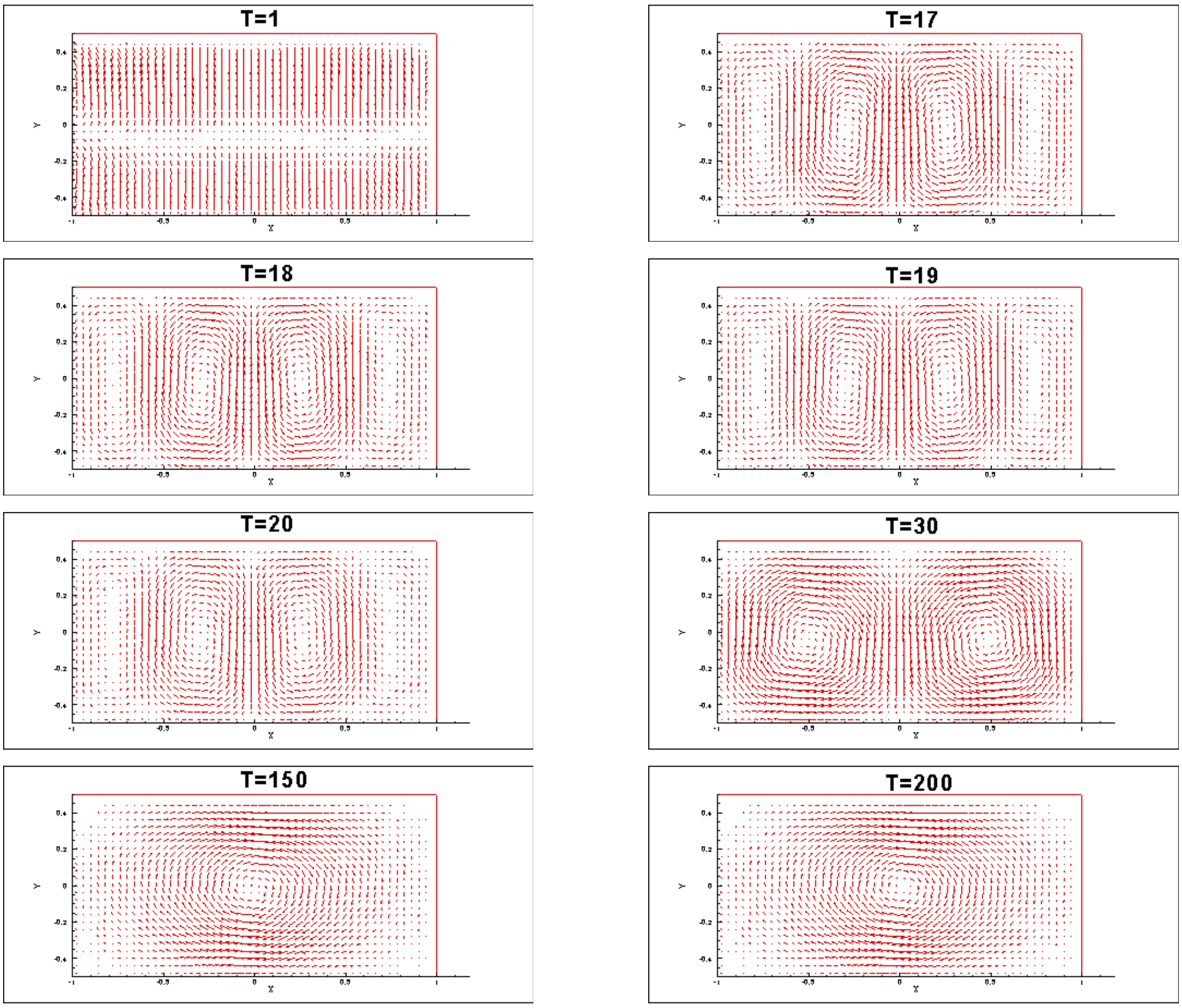

FIG. 7. The dynamical evolution of the velocity field $\boldsymbol{u}$ under external the magnetic field with the time step $\delta t = 0.01$. Snapshots of the numerical approximation are taken at $t = 1, 17, 18, 19, 20, 30, 150, 200$.

**Acknowledgments.** R. Chen is partially supported by the China Postdoctoral Science Foundation grant No. 2016M591122. X. Yang is partially supported by National Science Foundation under grant numbers DMS-1418898 and DMS-1720212. The work of H. Zhang is partially supported by NSFC-11471046, NSFC-11571045 and the Fundamental Research Funds for the Central Universities.